\documentclass[11pt]{article}
\usepackage{anysize}
\marginsize{2.0cm}{2.0cm}{2.0 cm}{2.0 cm}
\usepackage{setspace}

\usepackage{latexsym,amsfonts,amsmath}
\usepackage{graphicx}
\usepackage{color}
\usepackage{epstopdf}
\usepackage[titletoc]{appendix}

\def\liminf{\mathop{\underline{\lim}}}

\newtheorem{condition}{Condition}[section]{\bfseries}{\itshape}

{\bfseries}{\itshape}

\newtheorem{theorem}{Theorem}[section]{\bfseries}{\itshape}

\newtheorem{corollary}{Corollary}[section]{\bfseries}{\itshape}

\newtheorem{proposition}{Proposition}[section]{\bfseries}{\itshape}

\newtheorem{example}{Example}[section]{\bfseries}{\itshape}

\newtheorem{lemma}{Lemma}[section]{\bfseries}{\itshape}

\newtheorem{remark}{Remark}[section]{\bfseries}{\itshape}

\newtheorem{definition}{Definition}[section]{\bfseries}{\itshape}

{\bfseries}{\itshape}

\begin{document}
\title{On gradual-impulse control of continuous-time Markov decision processes with exponential utility}

\author{Xin Guo\thanks{Department of Mathematical Sciences, University of
Liverpool, Liverpool, U.K.. E-mail: x.guo21@liv.ac.uk.}, Aiko Kurushima\thanks{Department of Economics, Sophia University, Tokyo, Japan. Email: kurushima@sophia.ac.jp.}, Alexey Piunovskiy\thanks{Department of Mathematical Sciences, University of
Liverpool, Liverpool, U.K.. E-mail: piunov@liv.ac.uk.}~ and Yi
Zhang \thanks{Department of Mathematical Sciences, University of
Liverpool, Liverpool, L69 7ZL, U.K.. E-mail: yi.zhang@liv.ac.uk.}}
\date{}
\maketitle

\par\noindent{\bf Abstract:} In this paper, we consider the gradual-impulse control problem of continuous-time Markov decision processes, where the system performance is measured by the expectation of the exponential utility of the total cost. We prove, under very general conditions on the system primitives, the existence of a deterministic stationary optimal policy out of a more general class of policies. Policies that we consider allow multiple simultaneous impulses, randomized selection of impulses with random effects, relaxed gradual controls, and accumulation of jumps. After characterizing the value function using the optimality equation, we reduce the continuous-time gradual-impulse control problem to an equivalent simple discrete-time Markov decision process, whose action space is the union of the sets of gradual and impulsive actions.
\bigskip

\par\noindent {\bf Keywords:} Continuous-time Markov decision
processes. Exponential utility. Impulse-gradual control. Risk-sensitive criterion.  Optimality equation.
\bigskip

\par\noindent
{\bf AMS 2000 subject classification:} Primary 90C40,  Secondary
60J75
\section{Introduction}\label{JapanSec01}
This paper considers the gradual-impulse control problem for continuous-time Markov decision processes (CTMDPs) with the performance to be minimized being the exponential utility of the total cost. In this model, the decision maker can control the process gradually via its local characteristics (transition rate), and also has the option of affecting impulsively the state of the process. {The system dynamics is depicted in Figure \ref{Japanfig:01} below.} %Hence, the present paper is an extension of \cite{Zhang:2017}, where only gradual control was allowed.

There is no lack of situations, where an action can affect the state of the controlled process instantaneously. For example, in a Susceptible-Infected-Recovered (SIR) epidemic model, the controller elaborates the immunization policy, affecting the transition rate from the susceptibles to the infectives, as well as the isolation policy, which reduces instantaneously the number of infectives. {Let us formulate another simple example, which contains some features motivating the present paper.
\begin{example}\label{JapanRevisionExample001}
A rat (or intruder) may invade the kitchen. For each time unit it remains alive in the ``kitchen'', a constant cost of $l\ge 0$ is incurred. The rat spends an exponentially distributed amount of time with mean $\frac{1}{\mu}>0$ in the kitchen, and then goes outside and settles down in another house (and thus never returns). When the rat is in the kitchen, the housekeeper (defender) can decide to shoot at it, with a chance of hitting and killing the rat being $p\in(0,1)$. If the rat dodged, it remains in the kitchen. Each bullet costs $C>0$. Assume that the successive shootings are independent.
\end{example}
Let us mention some features in the above example. ``Shoot'' is an impulse. The location of the rat is the state. The effect of an impulse on the post-impulse state is random, as the shooting may be dodged. It is costly for each time unit the rat is present in the kitchen. Suppose the cost of impulse is relatively low. It can happen that after one impulse, if the rat is still alive and in the kitchen, then it is reasonable to immediately shoot again. This means, one should allow multiple impulses at a single time moment in this problem. We will return to this problem in Example \ref{JapanRevisionExample005} below, which demonstrates the situations when applying only one impulse is insufficient for optimality.}

Most previous works on gradual-impulse control do not allow multiple simultaneous impulses at a single time moment, see  \cite{Costa:1989,Costa:2000,Davis:1993,Hordijk:1984,Miller:2018,Palczewski:2017,vanderDuynSchouten:1983}. {Extra conditions are imposed therein to guarantee there is no need to apply more than one impulse at a single time moment. Example \ref{JapanRevisionExample001} described a situation when those conditions are not satisfied. There is convenience if only one impulse is allowed at a given time moment, because at each time moment, there is only one state, so that one can construct the process under control in the (original) state space of the gradual-impulse control problem. When there is only gradual control, it is convenient to construct the CTMDP using a marked point process with the mark space being the same as the state space of the original control problem.
%It is convenient to construct a CTMDP with gradual control using an underlying marked point process, whose mark space is the same as the state space of the process.
If multiple impulses were applied in a sequence at a single time moment, then there would be multiple states associated with the single time moment. If one wishes to construct the problem using a marked point process, then the mark space must be enlarged, so that a sequence of impulses applied at the single time moment and the post-impulse states are merged as a single ``mark'', which will be called an intervention.  Necessarily this leads to a complicated marked point process with extra notations, with each mark corresponding to a sample path of a discrete-time Markov decision process (DTMDP). This idea was employed and implemented in \cite{DufourPiunovskiy:2016}.}

Another way of constructing rigorously a gradual-impulse control problem of CTMDPs admitting multiple simultaneous impulses comes from \cite{Yushkevich:1988}. The idea is to keep the original state space, but to enlarge the time $t\in[0,\infty)$ to $(n,t)$ with the first coordinate, roughly speaking, counting the number of impulses applied at the time $t$. Consequently, several concepts about stochastic processes needed be extended. % Another way of mitigating the construction problem arising from multiple simultaneous impulses is to merge the sequence of impulses applied at the single time moment and the post-impulse states as a single ``state'', which will be called intervention. After that, the gradual-impulse control problem can be described using the standard theory of marked point processes. This idea is successfully applied in Dufour and Piunovskiy \cite{DufourPiunovskiy:2016} and its extention \cite{DufourHoriguchiPiunovskiy:2016}.

In the present work, we follow the construction of \cite{DufourPiunovskiy:2016} but with  more general control policies. Compared to the previous literature on impulse or gradual-impulse control problems of CTMDPs, to the best of our knowledge, we consider the most general setup: the policy allows to make relaxed gradual controls and randomized impulsive controls with randomized consequences, multiple simultaneous impulses are allowed, and accumulation of jumps of the process is not excluded. Another difference is that we consider the gradual-impulse control problem of CTMDPs with the system performance measure being the expectation of the exponential utility of the total cost to be minimized. %Problems with this performance measure are also called with multiplicative cost or risk-sensitive, as compared to the linear utility case, which is called risk-neutral. One of the pioneering works on risk-sensitive control appeared in 1970s, see \cite{Howard:1972}, where the justification of use of the term ``risk-sensitive'' was provided, and there have been reviving interest in it in the recent two decades. For DTMDP problems, see e.g., \cite{Cavazos:2000,Cavazos:2000b,Di Masi:1999,Jaskiewicz:2008,Patek:2001} and \cite{BauerleRieder:2014,Haskell:2015}: the latter references consider a more general utility function.
For risk-sensitive CTMDPs with gradual control only, see e.g., \cite{Ghosh:2014,XinGuo:2018,Kumar:2013,Piunovski:1985,Wei:2016,Zhang:2017}. In close relation to the present paper, the risk-sensitive optimal stopping problem of a continuous-time Markov chain was recently considered in \cite{BauerlePopp:2018}, which is a special impulse control problem but with a more general utility function. There seems to be limited literature on risk-sensitive control of CTMDPs with both gradual and impulse actions.

The main optimality results of this paper lie in the following. We characterize the value function of the gradual-impulse control problem for CTMDPs with exponential utility in terms of the optimality equation, and show the existence of deterministic stationary optimal policies, under quite general and natural conditions compared to the literature.
For example, the growth on the gradual cost rates and impulse cost functions, as well as the transition rate can be quite general. In comparison, only bounded transition and cost rates were allowed in \cite{DufourPiunovskiy:2016}, which deals with a discounted problem with linear utility. The boundedness conditions were to guarantee that the Dynkin formula is applicable to functions of interest therein, which is important for the argument therein. %In addition, there are other interesting observations that can be drawn from the investigations.

%on gradual-impulse control problems of CTMDPs or piecewise deterministic processes, which often assume the boundedness of the transition or cost rates, see e.g., \cite{Costa:1989,Costa:2000,DufourPiunovskiy:2016,Hordijk:1984,Plum:1991}, because they either follow a different method (such as the infinitesimal method), or aim at obtaining different results (e.g., to show the value function as the unique solution to the quasi-variational inequalities). All these works were about problems with a risk-neutral criterion.

{The method of investigations in the present paper is different from \cite{DufourPiunovskiy:2016}, but is closer to \cite{Zhang:2017}, which studies a similar problem for CTMDPs but with gradual control only. Although both the present paper and \cite{Zhang:2017} follow the same idea of reducing the original problem to a DTMDP, the implementation for the gradual-impulse control problem becomes more involving. If only gradual control is allowed, this DTMDP is obtained by inspecting the continuous-time process after each jump moment, and the action is a decision rule, specifying the selection of actions between two consecutive jumps. The action space is thus the space of measurable mappings. If both gradual and impulsive controls are allowed, the DTMDP is obtained by inspecting the system dynamics after either a natural jump or an impulse, and the action is a triplet, including the time until the next time of applying an impulse (if no natural jump occurred then), the next impulse itself, and the decision rule for the selection of gradual controls. Compared to \cite{Zhang:2017}, the connection  between a strategy in the induced DTMDP and a policy in the gradual-impulse control, which is at the core of the justification of the reduction,  becomes more delicate, see Subsection \ref{JapanRevisionSubsection001} below.
}

{Apart from being with a complicated action space, the induced DTMDP model is not so convenient. For example, it is not a semicontinuous model even if the system primitives of the gradual-impulse control problem satisfy the compactness-continuity conditions, see Examples \ref{JapanRevisionExample002} and \ref{JapanExample01}. Consequently, the existence of an optimal policy does not follow automatically from the reduction to the DTMDP. Therefore, the second step in the investigation is to reduce further the DTMDP model to yet another one, which is a semicontinuous model, and with a simple action space (the union of the set of gradual actions and impulses). This second reduction is done based on the investigation of the optimality equation of the DTMDP obtained from the first reduction.
}

 The rest of the paper is organized as follows. We present the rigorous construction of the controlled process and problem statement in Section \ref{JapanSecModel}. Section \ref{JapanSecMain} consists of the main optimality results, whose proof is postponed to Section \ref{JapanSecProof}. The argument is based on the connection with a DTMDP model, which is introduced in Section \ref{JapanSecDTMDPhat}. To improve the readability, we  summarize the relevant notions and facts about DTMDPs in the Appendix.

% see e.g., \cite{Costa:1989,Costa:2000,DufourPiunovskiy:2016,DufourHoriguchiPiunovskiy:2016,Hordijk:1984,Gatarek:1992,Miller:2018,Plum:1991,vanderDuynSchouten:1983,Yushkevich:1983,Yushkevich:1988},
% Compared to the literature on impulse or gradual-impulse control problems of CTMDPs or PDPs, bounded transition and cost rates were often imposed, see \cite{Costa:2000,DufourPiunovskiy:2016,DufourHoriguchiPiunovskiy:2016,Hordijk:1984,Gatarek:1992}, because either a different method was followed or the objective was different.

%In \cite{Hordijk:1984,vanderDuynSchouten:1983}, the process evolves deterministically between two natural jumps.

% The more general class of processes than this is the piecewise deterministic processes (PDPs), see. In a PDP, apart from the natural jumps, there is also a possibility of deterministic jumps, which take place when the process hits a boundary of the state space.

% In comparison,where, like in most early literature, the authors described the problem formulation only formally because they allowed multiple impulses at a single time moment.

%The process under control in the present paper is a Markov pure jump process.

\bigskip

\par\noindent\textbf{Notations and conventions.} In what follows, ${\cal{B}}(X)$ is
the Borel $\sigma$-algebra of the topological space $X,$ $I$ stands for the indicator function, and $\delta_{x}(\cdot)$
is the Dirac measure concentrated on the singleton $\{x\},$ assumed to be measurable. A measure is $\sigma$-additive and $[0,\infty]$-valued.
Here and below, unless stated otherwise, the term of
measurability is always understood in the Borel sense. Throughout
this paper, we adopt the conventions of
\begin{eqnarray*}
\frac{0}{0}:=0,~0\cdot\infty:=0,~\frac{1}{0}:=+\infty,~\infty-\infty:=\infty.
\end{eqnarray*}
For each function $f$ on $X$, let $||f||:=\sup_{x\in X}|f(x)|.$

\section{Model description and problem statement}\label{JapanSecModel}
\subsection{System primitives of the gradual-impulse control problem}\label{JapanSubsec01}
We describe the primitives of the model as follows. The state space is $\textbf{X}$, the space of gradual controls is $\textbf{A}^G$, and the space of impulsive controls is $\textbf{A}^I$. It is assumed that $\textbf{X}$, $\textbf{A}^G$ and $\textbf{A}^I$ are all Borel spaces, endowed with their Borel $\sigma$-algebras ${\cal B}(\textbf{X}),$ ${\cal B}(\textbf{A}^G)$ and ${\cal B}(\textbf{A}^I)$, respectively. The transition rate, on which the gradual control acts, is given by $q(dy|x,a)$, which is a signed kernel from $\textbf{X}\times\textbf{A}^G$, endowed with its Borel $\sigma$-algebra, to ${\cal B}(\textbf{X}),$ satisfying the following conditions: $q(\Gamma|x,a)\in[0,\infty)$ for each $\Gamma\in{\cal B}(\textbf{X}), x\notin \Gamma;$
\begin{eqnarray*}
&&q(\textbf{X}|x,a)=0,~x\in \textbf{X},~a\in\textbf{A}^G;~\bar{q}_x:=\sup_{a\in \textbf{A}^G}q_x(a)<\infty,~x\in\textbf{X},
\end{eqnarray*}
where $q_x(a):=-q(\{x\}|x,a)$ for each $(x,a)\in\textbf{X}\times\textbf{A}^G.$ For notational convenience, we introduce
\begin{eqnarray*}
\tilde{q}(dy|x,a):=q(dy\setminus\{x\}|x,a),~\forall~x\in\textbf{X},~a\in\textbf{A}^G.
\end{eqnarray*}
If the current state is $x\in\textbf{X}$, and an impulsive control $b\in\textbf{A}^I$ is applied, then the state immediately following this impulse obeys the distribution given by $Q(dy|x,b)$, which is a stochastic kernel from $\textbf{X}\times\textbf{A}^I$ to ${\cal B}(\textbf{X}).$ Finally, given the current state $x\in\textbf{X}$, the cost rate of applying a gradual control $a\in \textbf{A}^G$ is $c^G(x,a)$ and the cost of applying an impulsive control $b\in\textbf{A}^I$ is $c^I(x,b,y)$, where $c^G$ and $c^I$ are $[0,\infty)$-valued measurable functions on $\textbf{X}\times\textbf{A}^G$ and $\textbf{X}\times\textbf{A}^I\times\textbf{X}$, respectively.

Throughout this paper, we assume that $\textbf{A}^G$ and $\textbf{A}^I$ are compact Borel spaces. It is without loss of generality to assume $\textbf{A}^G$ and $\textbf{A}^I$ as two disjoint compact subsets of a Borel space $\tilde{\textbf{A}}$, for otherwise, one can consider $\textbf{A}^G\times\{G\}$ instead of $\textbf{A}^G$ and $\textbf{A}^I\times\{I\}$ instead of $\textbf{A}^I$ and $\tilde{\textbf{A}}=\textbf{A}^G\times\{G\}\bigcup \textbf{A}^I\times\{I\}$.
Furthermore, we assume that
\begin{eqnarray}\label{JapanBounded01}
\sup_{a\in\textbf{A}^G}c^G(x,a)<\infty,~\forall x\in\textbf{X}.
\end{eqnarray}
In what follows, we will not make specific reference to this assumption.

{The system dynamics in the concerned gradual-impulse control problem can be described as follows. In absence of impulses, the system is just a controlled Markov pure jump process in the state space $\textbf{X}$, where the (gradual) control, selected from $\textbf{A}^G$, acts on the local characteristics of the process, leading to natural jumps. This is conveniently described as a marked point process, which consists of the pairs of subsequent jump moments and the the post-jump states (marks). The mark space is thus $\textbf{X}$. We would still describe the system in the concerned gradual-impulse control problem using a marked point process. However, when the decision maker is allowed to apply a finite or countably infinite sequence of impulses from $\textbf{A}^I$ at a single time moment, and each impulse results in a post-impulse state, there would be a sequence of states in $\textbf{X}$ at a single time moment. Moreover, the order of the impulses and their resulting states is also relevant. Therefore, the marked point process we use now is in an enlarged mark space. More precisely, each mark contains a sequence of impulses applied at the same time moment, the state before the impulses are applied, and all the states resulted  by these impulses. Each jump moment is either triggered by an impulse (or a sequence of impulses), or by a natural jump. A mark in this marked point process is referred to as an intervention. This term is naturally understandable when the mark consists of impulses. Having said so, we will also allow that an ``intervention'' does not contain any impulse or say an empty sequence of impulses. This appears when the decision maker chooses not to apply any impulse immediately after a natural jump.}
In the rest of this section, following the method of \cite{DufourPiunovskiy:2016}, we will elaborate this idea and describe rigorously the concerned continuous-time gradual-impulse control problem. {To this end, we will firstly state the precise definition of an intervention in the next subsection.} %and view the sequence of impulses applied at the single time moment and the post-impulse state as a single ``state'', which will be called intervention.

%As mentioned in the Introduction, another way of rigorously defining the controlled process under multiple impulses at a single time moment was given in \cite{Yushkevich:1988}. Here we choose to follow the construction in \cite{DufourPiunovskiy:2016} for its simplicity and generality (we consider randomized and relaxed policy, whereas \cite{Yushkevich:1988} considered only deterministic policies).

\subsection{Definition and interpretation of an intervention}
{
At the beginning of an intervention, the decision marker chooses whether to apply an impulse, and which one to apply. If the current state is $x\in\textbf{X}$, and after an impulse $b\in\textbf{A}^I$ is chosen, the new state say $y\in\textbf{X}$ is instantaneously realized, following the distribution $Q(dy|x,b)$. Then based on $x,b,y$, the decision maker will choose the next impulse, if any at all, and so on. To be consistent, a cemetery point $\Delta\notin \textbf{A}^I,\textbf{X}$ is artificially fixed, which is chosen when the decision maker decides not to apply any more impulse at the current time, and it leads to the post-impulse state, also denoted as $\Delta,$ which is absorbing, i.e., $Q(\Delta|\Delta,\Delta)\equiv 1.$ Therefore, an intervention is a sequential decision process. More precisely, an intervention can be regarded as a trajectory or sample path of the following DTMDP, which we refer to as the ``intervention'' DTMDP model, to distinguish it from several other DTMDP models to appear subsequently.
\begin{definition}
The intervention DTMDP model is specified by the following tuple $\{\mathbf{X}_{\Delta},\mathbf{A}^{I}_{\Delta},Q\}$, which are defined in terms of the primitives of the gradual-impulse control problem given in Subsection \ref{JapanSubsec01}.
\begin{itemize}
\item The state space is $\mathbf{X}_{\Delta}:=\mathbf{X}\bigcup\{\Delta\}$, where  $\Delta$ is a cemetery point not belonging to $\textbf{X}$ or $\textbf{A}^I.$
\item The action space is  $\mathbf{A}^{I}_{\Delta}:= \mathbf{A}^{I}\bigcup\{\Delta\}.$
\item The one-step transition probability from $\textbf{X}_\Delta\times\textbf{A}^I_\Delta$ to ${\cal B}(\textbf{X}_\Delta)$ is $Q(dy|x,b)$, where we have accepted that $Q(\{\Delta\}|x,b):=1$ if $x=\Delta$ or $b=\Delta$.
\end{itemize}
\end{definition}
}

Let the initial distribution in the intervention DTMDP be always concentrated on $\textbf{X}.$ Then its canonical sample space is
\begin{eqnarray*}
\mathbf{Y}:=\left(\bigcup_{k=0}^{\infty}\mathbf{Y}_k\right) \bigcup(\textbf{X}\times\textbf{A}^I)^\infty,
\end{eqnarray*}
where for each $\infty>k\ge 1$
\begin{eqnarray*}
\mathbf{Y}_k:=(\mathbf{X}\times \mathbf{A}^I)^k\times(\mathbf{X}\times\{\Delta\})\times(\{\Delta\}\times\{\Delta\})^\infty,
\end{eqnarray*}
and $\textbf{Y}_0:=(\mathbf{X}\times\{\Delta\})\times(\{\Delta\}\times\{\Delta\})^\infty.$ Here, if $y\in \textbf{Y}_k$, $\infty>k\ge 0$, then there are $k$ impulses applied in the intervention $y$. %Similarly, if $y\in (\textbf{X}\times\textbf{A}^I)^\infty,$ then there are infinitely many impulses applied in the intervention $y.$

%An intervention is a sample path of the so-called intervention DTMDP with the following primitives:

Now we give the following definition.
{\begin{definition}
An intervention is an element of $\textbf{Y}.$
\end{definition}
In other words, $\textbf{Y}$ defined above is the space of all interventions. It will be the mark space of the marked point process $\{(T_n,Y_n)\}$ introduced in the next subsection.}

{With the notations introduced above, we now reiterate, more rigorously compared to the one in the beginning of this subsection, the interpretation of an intervention as follows}. Given the current state $x\in \textbf{X}$, if the controller decides to use $\Delta$, then it means, no more  impulse is used at this moment, and the intervention DTMDP is absorbed at $\Delta$; if the controller decides to use an impulse $b\in\textbf{A}^I$, then the post-impulse state follows the distribution $Q(dy|x,b)$. At the next post-impulse state $y$, if $y=\Delta,$ then the only decision is $\Delta;$ if $y\ne \Delta,$ then the controller either decides to use no impulse, leading to the next post-impulse state $\Delta$, or to use impulse $b'$, leading to the next post-impulse state, which follows the distribution given by $Q(\cdot|y,b')$, and so on. In other words, an intervention consists of a state and a finite or countable sequence of pairs of impulsive actions and the associated post-impulse states. In particular,
no impulse is applied in an intervention if the intervention belongs to $\mathbf{Y}_{0}$, {see Figure \ref{Japanfig:01} and its caption for an example.} Let \begin{eqnarray*}
\mathbf{Y}^{*}:=\textbf{Y}\setminus \textbf{Y}_0=\left(\bigcup_{k=1}^{\infty}\mathbf{Y}_k\right) \bigcup (\textbf{X}\times\textbf{A}^I)^\infty
\end{eqnarray*}
be the set of interventions, where some impulses are applied.

{In an intervention, locally, the selection of impulses (including the ``pseudo'' impulse $\Delta$) from $\textbf{A}^I_\Delta$ is governed by a strategy in the intervention DTMDP model. This adverb ``locally'' is understood in comparison with the definition of a policy for the gradual-impulse control problem, as  given in Definition \ref{JapanDef01} below, which governs the selection of impulsive controls as well as gradual controls, and is thus  ``global''. }
Let $\Xi$ be the set of (possibly randomized and history-dependent) strategies $\sigma$ in the intervention DTMDP.
We refer the reader to the appendix or \cite{Hernandez-Lerma:1996,Piunovskiy:1997} for standard terminologies in the theory of DTMDPs. {The way how a strategy in the intervention DTMDP model is incorporated into a policy in Definition \ref{JapanDef01} below is through its strategic measure. We recall the definition of a strategic measure in a DTMDP model in Definition \ref{JapanRevisionDef01}.} Let $\beta^\sigma(\cdot|x)$ denote the corresponding strategic measure of a strategy $\sigma$ of the intervention DTMDP, given the initial state $x\in\textbf{X}$. By the Ionescu-Tulcea theorem, see e.g., Proposition C.10 in \cite{Hernandez-Lerma:1996}, the mapping $x\in\textbf{X}\rightarrow \beta^\sigma(\cdot|x)$ is measurable. Let
$
\mathcal{P}^\mathbf{Y}
$
be the collection of all such stochastic kernels generated by some strategy $\sigma\in\Xi$, and
\begin{eqnarray*}
\mathcal{P}^\mathbf{Y}(x):=\{\beta^\sigma(\cdot|x): \sigma\in\Xi\}
\end{eqnarray*}
for each state $x\in \mathbf{X}$. Let
\begin{eqnarray*}
\mathcal{P}^{\mathbf{Y}^\ast}:=\{\beta(\cdot|\cdot)\in \mathcal{P}^{\textbf{Y}}:~\beta(\textbf{Y}^\ast|x)=1,~\forall~x\in\textbf{X}\},
\end{eqnarray*}
and for each $x\in\textbf{X},$
\begin{eqnarray*}
\mathcal{P}^\mathbf{Y^\ast}(x):=\{\beta(\cdot|x): \beta(\cdot|\cdot)\in \mathcal{P}^{\mathbf{Y}},~\beta(\textbf{Y}^\ast|x)=1\}.
\end{eqnarray*}

\subsection{Construction of the controlled processes}
{Let us now describe the promised marked point process $\{(T_{n},Y_{n})\}_{n=1}^\infty$ for the system dynamics of the concerned gradual-impulse control problem, where the mark space is the space of interventions.} Then the continuous-time controlled process $\{\xi_t\}_{t\ge 0}$ is defined based on the marked point process.

Let
\begin{eqnarray*}
&&\mathbf{Y}_\Delta:=\mathbf{Y}\bigcup\{\Delta\},\\
&& \Omega_0:=\mathbf{Y}\times(\{0\}\times \textbf{Y})\times(\{\infty\}\times\{\Delta\})^\infty,\\
&& \Omega_{n}:=\mathbf{Y}\times(\{0\}\times \textbf{Y})\times ((0,\infty)\times \mathbf{Y})^n\times(\{\infty\}\times\{\Delta\})^\infty,\forall~n=1,2,\dots.
\end{eqnarray*}
The canonical space $\Omega$ is defined as
\begin{eqnarray*}
\Omega:=\left(\bigcup_{n=0}^\infty \Omega_{n}\right)\bigcup \big( \mathbf{Y}\times((0,\infty)\times \mathbf{Y})^\infty \big)
\end{eqnarray*}
and is endowed with its Borel $\sigma$-algebra denoted by $\mathcal{F}$. The following generic notation of a point in $\Omega$ will be in use:
$
\omega=(y_0,\theta_1,y_1,\theta_2,y_2,\ldots).
$
Below, unless stated otherwise, $x_0\in\textbf{X}$ will be a fixed notation as the initial state of the  gradual-impulse control problem.
Then we put \begin{eqnarray}\label{JapanRevision01}
y_0:=(x_0,\Delta,\Delta,\dots),~\theta_1\equiv 0.
\end{eqnarray}
The sequence of $\{\theta_n\}_{n=1}^\infty$ represents the sojourn times between consecutive interventions. Here $\theta_1=0$ corresponds to that we allow the possibility of applying impulsive control at the initial time moment, c.f. (\ref{Japan0001}) below.

For each $n=0,1,\dots$, let \begin{eqnarray*}
h_{n}:=(y_0,\theta_1,y_1,\theta_2,y_2,\dots \theta_{n},y_{n})=(y_0,0,y_1,\theta_2,y_2,\dots \theta_{n},y_{n}),
\end{eqnarray*}
where the second equality holds because $\theta_1\equiv 0,$ see (\ref{JapanRevision01}).
 The collection of all such fragmental histories $h_n$ is denoted by $\mathbf{H}_{n}$. Let us introduce the coordinate mappings:
\begin{eqnarray*}
Y_{n}(\omega)=y_{n},~\forall~n\ge 0;~\Theta_{n}(\omega)=\theta_{n},~\forall~n\ge 1.
\end{eqnarray*}
The sequence $\{T_{n}\}_{n=1}^{\infty}$ of $[0,\infty]$-valued mappings is defined on $\Omega$ by
$
T_{n}(\omega):=\sum_{i=1}^n\Theta_i(\omega)=\sum_{i=1}^n\theta_i
$
and
$
T_\infty(\omega):=\lim_{n\to\infty}T_{n}(\omega).
$
Let $
H_{n}:=(Y_0,\Theta_1,Y_1,\dots,\Theta_{n},Y_{n}).
$
Finally, we define the controlled process $\big\{\xi_t\big\}_{t\in [0,\infty)}$:
\begin{eqnarray*}\xi_t(\omega)=\left\{
\begin{array}{ll}
Y_{n}(\omega), & \mbox{ if } T_{n}\le t<T_{n+1} \mbox{ for } n\ge 1; \\
\Delta, & \mbox{ if } T_\infty\le t,
\end{array}\right..
\end{eqnarray*}
%and $\xi_{0}(\omega)=Y_0=y_0$ with $y_0=(x_0,\Delta,\Delta,\dots)$.

It is convenient to introduce the random measure $\mu$ of the marked point process $\{(T_{n},Y_{n})\}_{n=1}^\infty$ on $(0,\infty)\times \mathbf{Y}$:
  \begin{eqnarray*}
  \mu(dt\times dy)=\sum_{n\ge 2}I_{\{T_{n}<\infty\}}\delta_{(T_{n},Y_{n})}(dt\times dy),
  \end{eqnarray*}
where the dependence on $\omega$ is not explicitly indicated.
Let
$\mathcal{F}_t:=\sigma\{H_1\}\vee\sigma\{\mu((0,s]\times B):~s\le t,B\in\mathcal{B}(\mathbf{Y})\}$
for $t\in[0,\infty)$.

\begin{figure}
  \includegraphics[width=1.0\textwidth, angle = 0]{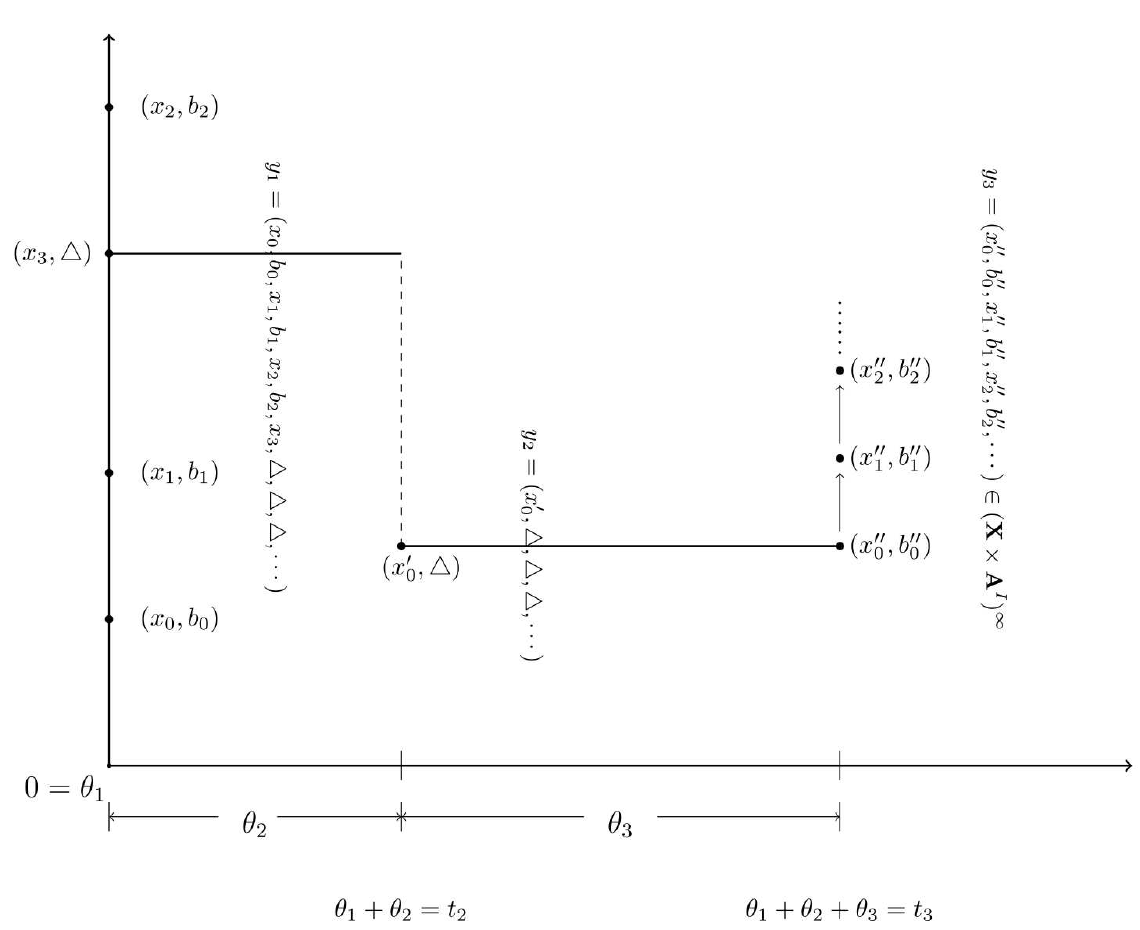}
  \caption{\small {Illustration of the system dynamics in the gradual-impulse control problem, and how the policy acts on the system dynamics. Here $\textbf{X}=[0,\infty).$ The second coordinate indicates the impulse (including the ``pseudo'' impulse $\Delta$) used at that state, which is recorded in the first coordinate. At the initial time $t=\theta_1\equiv 0$, three impulses are applied in turn. The first jump in the indicated sample path of the marked point process $\{(T_n,Y_n)\}_{n=1}^\infty$ takes place at $t_2=\theta_2.$ It is triggered by a natural jump because $x_0'\ne x_3.$ Along the displayed sample path, the system state remains to be $x_3$ before the first jump of the marked point process. The second jump of the marked point process is triggered by a planned (or say active) impulse, because $x_0''=x_0'$. Infinitely many impulses are applied at $t_3=t_2+\theta_3$, so that the process is ``killed'' after the ``infinitely many'' impulses at $t_3,$ i.e., $\omega=(y_0,0,y_1,\theta_2,y_2,\theta_3,y_3,\infty,\Delta,\infty,\Delta,\dots).$ Note also that, under the policy $u=\{u_{n}\}_{n=0}^\infty$ in Definition \ref{JapanDef01}, $y_1\in \textbf{Y}_3$ is a realization from the distribution $u_0(\cdot|x_0)$, $\bar{x}(y_1)=x_3;$ $y_2\in \textbf{Y}_0$ is a realization from the distribution $\Gamma^0_1(\cdot|h_1,\theta_2,x_0')$ as the jump at $t_2$ is triggered by a natural jump, $\bar{x}(y_2)=x_0'$; and $y_3\in (\textbf{X}\times\textbf{A}^I)^\infty$ is a realization from the distribution $\Gamma_2^1(\cdot|h_2)$ as the jump at $t_3$ is not triggered by a natural jump, $\bar{x}(y_3)=\Delta.$} }
  \label{Japanfig:01}
\end{figure}

We will use the following notation in the next definition. For each $y=(x_0,b_0,x_1,b_1,\dots)\in \mathbf{Y}$,
\begin{eqnarray*}
\bar x(y):=x_{k}
\end{eqnarray*}
if $\infty>k=0,1,\dots$ is the unique integer such that $y\in \mathbf{Y}_k$ (if $k\ge 1,$ then $\bar{x}(y)$ is the state after the last impulse in the intervention $y$); if such an integer $k$ does not exist, then $y\in(\textbf{X}\times\textbf{A}^I)^\infty$ and \begin{eqnarray*}
\bar{x}(y):=\Delta.
\end{eqnarray*}
That previous equality corresponds to that we kill the process after an infinite number of impulses was applied at a single time moment. {An example of a trajectory of the system dynamics in the gradual impulse control problem is displayed in Figure \ref{Japanfig:01}.
}

\begin{definition}\label{JapanDef01}
A policy is a sequence $u=\{u_{n}\}_{n=0}^\infty$ such that
$u_{0}\in \mathcal{P}^{\mathbf{Y}}$ and, for each $n=1,2,\dots$,
\begin{eqnarray*}
u_{n}=\left( \Phi_{n},\Pi_{n},\Gamma_{n}^0,\Gamma_n^1 \right),
\end{eqnarray*}
where
$\Phi_{n}$ is a stochastic kernel on $(0,\infty]$ given $\mathbf{H}_{n}$,
$\Pi_{n}$ is a stochastic kernel on $\mathbf{A}^{G}$ given $ \mathbf{H}_{n}\times (0,\infty)$ such that $\Phi_n(\{\infty\}|h_n)=1$ if $y_n\in (\textbf{X}\times\textbf{A}^I)^\infty$,
$\Gamma^0_{n}$ is a stochastic kernel on $\mathbf{Y}$ given $ \mathbf{H}_{n}\times (0,\infty)\times \mathbf{X}$ satisfying $\Gamma_n^0(\cdot|h_n,t,x)\in {\cal P}^{\textbf{Y}}(x)$ for each $h_n\in\textbf{H}_n$ and $x\in\textbf{X}$ and $t\in(0,\infty)$; and $\Gamma^1_{n}$ is a stochastic kernel on $\mathbf{Y}$ given $ \mathbf{H}_{n}$ satisfying $\Gamma^1_n(\cdot|h_n)\in {\cal P}^{\textbf{Y}^\ast}(\bar{x}(y_n))$ for each $h_n\in \textbf{H}_n$.
(The above conditions apply when $y_n\ne \Delta$; otherwise, all the values of $\Phi_n(\cdot| h_n)$, $\Pi_n(\cdot|h_n,t)$, $\Gamma_n^0(\cdot|h_n,t,\cdot)$ are immaterial and may be put arbitrarily. )
\end{definition}
The set of policies is denoted by $\cal U$.

{Let us provide an interpretation of how a policy $u$ acts on the system dynamics.} Roughly speaking, an intervention is over as soon as the (possibly empty) sequence of simultaneous impulses is over. Given that the $n$th intervention is over, the kernel $\Phi_n$ specifies the conditional distribution of the planned time until the next impulse (or next sequence of impulses). The (conditional) distribution of the time until the next natural jump (if there were no interventions before it) is the non-stationary exponential distribution with rate $\int_{\textbf{A}^G}q_{\bar{x}(Y_n)}(a)\Pi_n(da|H_n,t)$. In other words, $\Pi_n$ is the relaxed gradual control. Given the $n$th intervention is over, the next intervention is triggered by either the next planned impulse or the next natural jump; in the former case, the new intervention has the distribution given by $\Gamma^1_n$, and in the latter case the new intervention has the distribution given by $\Gamma^0_n.$ This interpretation will be seen consistent with (\ref{JapanExplain01}) and (\ref{JapanExplain02}) below, {where one can see how a policy $u$ acts on the conditional law of the marked point process $\{(T_n,Y_n)\}_{n=1}^\infty.$ See also the caption of Figure \ref{Japanfig:01}.}

%Below, in most of the cases, the term of a control policy is associated with continuous-time impulse-gradual control problem. This is to distinguish it from the corresponding object in DTMDPs, which will be called strategy.

Suppose a policy $u=\{u_{n}\}_{n=0}^\infty $ is fixed. {Let us now present the conditional law of the marked point process $\{(T_n,Y_n)\}_{n=1}^\infty$ under the policy $u$, which determines the underlying probability measure ${\rm P}_{x_0}^u$ on $(\Omega,\cal F)$, where $x_0\in\textbf{X}$ is the fixed initial state of the system dynamic.} For brevity,
we introduce the following notations for each $n\ge 1$, $\Gamma\in \mathcal{B}(\mathbf{X})$ and $h_{n}=(y_0,\theta_1,y_1,\ldots,\theta_{n},y_{n})\in \mathbf{H}_{n}$:
\begin{eqnarray*}
&&\lambda_{n}^u(\Gamma|h_{n},t):  =  \int_{\mathbf{A}^{G}} \tilde{q}(\Gamma | \overline{x}(y_{n}),a) \Pi_{n}(da | h_{n},t),
\\
&&\Lambda_{n}^u(\Gamma|h_{n},t):  =  \int_{0}^t \lambda_{n}^u(\Gamma|h_{n},s) ds.
\end{eqnarray*}
where and below, we put $q_{\Delta}(a):=0$ for each $a\in\textbf{A}^G.$ Now, for each $n\ge 1$, we introduce the stochastic kernel $G_{n}^u$ on $(0,\infty]\times \mathbf{Y}_{\Delta}$ given $\mathbf{H}_{n}$ as follows. For each $h_{n}=(y_0,\theta_1,y_1,\ldots,\theta_{n},y_{n})\in \mathbf{H}_{n}$,
\begin{eqnarray}\label{JapanExplain01}
G_{n}^u(\{+\infty\}\times \{\Delta\} | h_{n}):= \delta_{y_{n}} (\{\Delta\}) + \delta_{y_{n}} (\mathbf{Y})
e^{-\Lambda_{n}^u(\mathbf{X}|h_{n},+\infty)}\Phi_{n}(\{+\infty\}|h_n),
\end{eqnarray} and
\begin{eqnarray}\label{JapanExplain02}
&&G_{n}^u(dt\times dy| h_{n}):=\delta_{y_{n}} (\mathbf{Y}) \left\{
  \Gamma_{n}^1(dy| h_{n}) e^{-\Lambda_{n}^u(\mathbf{X}|h_{n},t)} \Phi_{n}(dt | h_{n}) \right.\nonumber \\
&& \left.+  \int_{\mathbf{X}} \Phi_{n}([t,\infty] | h_{n}) \Gamma_{n}^0(dy| h_{n},t,x)  \lambda_{n}^u(dx|h_{n},t) e^{-\Lambda_{n}^u(\mathbf{X}|h_{n},t)} dt \right\}
\end{eqnarray}
on $(0,\infty)\times \textbf{Y}$. For each fixed initial state $x_{0}\in\mathbf{X}$, by the Ionescu-Tulcea theorem, see e.g., Proposition C.10 in \cite{Hernandez-Lerma:1996}, there exists a probability $ {\rm P}^{u}_{x_{0}}$ on $(\Omega,\mathcal{F})$ such that
the restriction of ${\rm P}^{u}_{x_{0}}$ to $(\Omega,\mathcal{F}_{0})$ is given by
\begin{eqnarray}\label{Japan0001}
{\rm P}^{u}_{x_{0}} \left( \left(\{y_{0}\}\times \{0\} \times \Gamma \times ((0,\infty]\times \mathbf{Y}_{\Delta})^{\infty} \right)\bigcap \Omega\right) & = &
u_{0}(\Gamma|x_{0})
\end{eqnarray}
for each $\Gamma\in \mathcal{B}(\mathbf{Y})$;  and for each $n\ge 1$, under ${\rm P}^{u}_{x_{0}}$,  the conditional distribution of $(Y_{n+1},\Theta_{n+1})$ given $\mathcal{F}_{T_{n}}:=\sigma(H_{n})$ is determined by
$G_{n}^u(\cdot | H_{n})$ and the conditional survival function of $\Theta_{n+1}$ given $\mathcal{F}_{T_{n}}$ under ${\rm P}^{u}_{x_{0}}$
is given by  $G_{n}^u([t,+\infty]\times \mathbf{Y}_{\infty}| H_{n})$.

The cost associated with an intervention  $y=(x_{0},b_{0},x_{1},b_{1},\ldots)\in \mathbf{Y}$ is given by
\begin{eqnarray*}
C^I(y):=\sum_{k=0}^\infty c^I(x_k,b_k,x_{k+1}).
\end{eqnarray*}
Here, recall that an intervention consists of the current state, the sequence of impulses applied in turn at the same time moment and the associated post-impulse states; and each impulse $b$ applied at state $x$ results in a cost $c^I(x,b,z)$ if it leads to the post-impulse state $z$. (We accept that $c^I(x,\Delta,\Delta):= 0$ for all $x\in\textbf{X}_\Delta$.) With this notation, we now introduce the performance measure considered in this paper:
{
\begin{eqnarray*}
\mathcal{V}(u,x)&:=& {\rm E}^{u}_{x} \left[ e^{\sum_{n=1}^\infty  \left(C^{I}(Y_n)+ \int_{T_n}^{T_{n+1}}  \int_{\mathbf{A}^{G}}  c^{G}(\bar{x}(\xi_{s}),a) \Pi_n(da |H_n,s-T_n) ds\right)} \right]
\end{eqnarray*}
for each $x\in\textbf{X}$ and policy $u\in {\cal U}$.  Here we recall that $T_1=\Theta_1\equiv 0$, see (\ref{JapanRevision01}). To illustrate more explicitly how the policy acts on the impulses,  consider the example of only one intervention and null gradual cost $c^G(x,a)\equiv 0$. Then we may write
\begin{eqnarray*}
{\rm E}^{u}_{x} \left[ e^{   C^{I}(Y_1)  } \right]=\int_{\textbf{Y}} u_0(dx_0\times db_0\times dx_1\times db_1\times\dots|x)e^{\sum_{k=0}^\infty c^I(x_k,b_k,x_{k+1}) }=\int_\textbf{Y} u_0(dy|x)e^{ C^I(y) }.
\end{eqnarray*}
More generally, one can compute ${\rm E}^{u}_{x} \left[ e^{   C^{I}(Y_{n+1})  } \right]={\rm E}^{u}_{x} \left[{\rm E}^{u}_{x} \left[ e^{   C^{I}(Y_{n+1})  }|H_n \right]\right]$, where ${\rm E}^{u}_{x} \left[ e^{   C^{I}(Y_{n+1})  }|H_n \right]$ can be written out as a similar integral to the case of $n=0$ using the conditional laws (\ref{JapanExplain01}) and (\ref{JapanExplain02}).
}

Let the value function ${\cal V}^\ast$ be denoted by
\begin{eqnarray*}
{\cal V}^\ast(x):=\inf_{u\in {\cal U}}{\cal V}(x,u)
 \end{eqnarray*}
 for each $x\in\textbf{X}$. A policy $u^\ast$ satisfying ${\cal V}(x,u^\ast)={\cal V}^\ast(x)$ for all $x\in\textbf{X}$ is called optimal for the gradual-impulse control problem:
\begin{eqnarray}\label{JapanProblem01}
\mbox{Minimize over $u\in{\cal U}:$~}{\cal V}(x,u).
\end{eqnarray}
In this paper, we will present conditions on the system primitives that guarantee the existence of an optimal policy in a simple form as defined as follows.
%We end this section with the following definition of a deterministic stationary policy.
\begin{definition}\label{JapanDef02}
A policy $u$ is called deterministic stationary if there exist some measurable mappings $(\varphi,\psi,f)$ on $\textbf{X}$, where $\varphi(x)\in\{0,\infty\}$ for each $x\in \textbf{X}$, $\psi$ and $f$ are $\textbf{A}^I$-valued and $\textbf{A}^G$-valued, such that
$\Phi_n(\{\infty\}|h_n)=1$, $\Pi_n(da|h_n,t)=\delta_{f(\bar{x}(y_n))}(da)$ for all $t\ge 0$, and $u_0(\cdot|x)=\Gamma^0_n(\cdot|h_n,t,x)=\beta^\pi(\cdot|x)$ for some deterministic stationary strategy $\pi$ in the intervention DTMDP model defined by $\pi(\{\Delta\}|x_0,b_0,x_1,b_1,\dots,x_n)=I\{\varphi(x_n)=\infty\}$, and $\pi(db|x_0,b_0,x_1,b_1,\dots,x_n)=I\{\varphi(x_n)=0\}\delta_{\psi(x_n)}(db).$
\end{definition}
In the above definition, $\Gamma^1_n$ was left arbitrary, because, under such a deterministic stationary policy, a new intervention is always triggered by a natural jump.

%We will also characterize the value function ${\cal V}^\ast$ using the optimality equation for a DTMDP problem, whose action space is the union of the sets of gradual and impulse actions. {To this end, we shall firstly relate the gradual-impulse control problem to a DTMDP model, through which we establish the optimality results for the original problem. This auxiliary DTMDP model is termed as the hat DTMDP model, and is introduced in the next section.}

\section{Optimality result}\label{JapanSecMain}
In this section, we present the main optimality results in this paper. In a nutshell, under quite general conditions on the system primitives of the gradual-impulse control problem (\ref{JapanProblem01}), we show that it can be solved via problem (\ref{ZyExponentialProblem2}) for a simple DTMDP model, which we refer to as the tilde DTMDP model. In this way, we show that the gradual-impulse control problem (\ref{JapanProblem01}) admits a deterministic stationary optimal policy. %This reduction is based on the characterization of the value function using optimality equations. We do this by analyzing the hat DTMDP problem (\ref{JapanProblem02}), which is in general more complicated than the tilde DTMDP problem.

In order to formulate the tilde DTMDP model, we impose the following condition.

\begin{condition}\label{JapanCon02}
There exists an $[1,\infty)$-valued continuous function $w$ on $\textbf{X}$ such that $c^G(x,a)+q_x(a)+1\le w(x)$ for each $(x,a)\in\textbf{X}\times\textbf{A}^G.$
\end{condition}
If $c^G$ is a continuous function, then the above condition is a consequence of Condition \ref{JapanCon01} and the Berge theorem, see Proposition 7.32 of \cite{Bertsekas:1978}. Several statements below do not need the bounding function $w$ in Condition \ref{JapanCon02} to be continuous. In this connection, we also mention that a Borel measurable function $w$ satisfying the inequality in Condition \ref{JapanCon02} always exists, see Lemma 1 of \cite{Feinberg:2017} and recall (\ref{JapanBounded01}).

{Recall that $\tilde{\textbf{A}}=\textbf{A}^I\bigcup\textbf{A}^G$ is the disjoint union of $\textbf{A}^G$ and $\textbf{A}^I$. We are now in position to define the tilde DTMDP model in terms of the system primitives of the gradual-impulse control problem (\ref{JapanProblem01}).
\begin{definition}
The tilde DTMDP model is specified by the following four-tuple $\{\textbf{X},\tilde{\textbf{A}},\tilde{Q},\tilde{l}\}$, where $\textbf{X}$ and $\tilde{\textbf{A}}$ are its state and action spaces, and its transition probability $\tilde{Q}$ on $\textbf{X}$ given $\textbf{X}\times\tilde{\textbf{A}}$  and cost function $\tilde{l}$ are defined by
\begin{eqnarray*}
\tilde{Q}(dy|x,a):= \frac{q(\Gamma|x,a)}{w(x)}+\delta_{x}(dy),~\tilde{l}(x,a,y):=\ln\frac{w(x)}{w(x)-c^G(x,a)}
 \end{eqnarray*} for all $a\in\textbf{A}^G$, \begin{eqnarray*}
 \tilde{Q}(dy|x,b):=Q(dy|x,b),~\tilde{l}(x,b,y):=c^I(x,b,y)
 \end{eqnarray*} for all $b\in\textbf{A}^I$.
\end{definition}
}

%As mentioned in the beginning of this section, we will establish the existence of a deterministic stationary optimal policy for the gradual-impulse control problem (\ref{JapanProblem01}) via the one of problem (\ref{ZyExponentialProblem2}) for the tilde DTMDP model.
For the solvability of problem (\ref{ZyExponentialProblem2}) for the tilde DTMDP model, we impose the following compactness-continuity condition.
\begin{condition}\label{JapanCon01}
The functions $c^I$ and $c^G$ are lower semicontinuous on $\textbf{X}\times\textbf{A}^I\times\textbf{X}$ and $\textbf{X}\times\textbf{A}^G$, respectively; and for each bounded continuous function $g$ on $\textbf{X}$, $\int_{\textbf{X}}g(y)Q(dy|x,b)$ and $\int_{\textbf{X}}g(y)\tilde{q}(dy|x,a)$ are continuous in $(x,b)\in \textbf{X}\times\textbf{A}^I$ and $(x,a)\in \textbf{X}\times\textbf{A}^G$, respectively. (Recall also that $\textbf{A}^G$ and $\textbf{A}^I$ are compact.)
\end{condition}
Under Conditions \ref{JapanCon02} and \ref{JapanCon01}, one can easily check that the tilde DTMDP model is semicontinuous, so that the value function $W^\ast$ for problem (\ref{ZyExponentialProblem2}) of the tilde DTMDP model is lower semicontinuous, and there exists an optimal deterministic stationary strategy for it, see Proposition \ref{GGZyExponentialProposition02}(e). We collect these observations in the next statement for future reference.

\begin{proposition}\label{JapanRevisionAAProp01}
Suppose Conditions \ref{JapanCon02} and \ref{JapanCon01} are satisfied. Then  the value function $W^\ast$ of problem (\ref{ZyExponentialProblem2}) for the tilde DTMDP model coincides is the minimal $[1,\infty]$-valued lower semicontinuous function satisfying
\begin{eqnarray}\label{JapanBellman05}
V(x)= \inf_{\tilde{a}\in\tilde{\textbf{A}}}\left\{\int_{\textbf{X}}e^{\tilde{l}(x,\tilde{a},y)}V(y)\tilde{Q}(dy|x,\tilde{a})\right\},~x\in\textbf{X},
\end{eqnarray}
and the above relation holds with equality being replaced by ``$\ge $'', too.
A pair of measurable mappings $(\psi^\ast,f^\ast)$ from $\textbf{X}$ to $\textbf{A}^I$ and $\textbf{A}^G$, respectively, % satisfies the relations in Theorem \ref{JapanTheorem01}(c) if and only if they are conserving in the tilde model, i.e., for each $x\in\textbf{X}$, there is some $\tilde{a}^\ast\in\tilde{\textbf{A}}$ such that
is a deterministic optimal stationary strategy for problem (\ref{ZyExponentialProblem2}) of the tilde DTMDP model if and only if
\begin{eqnarray}\label{JapanConserving01}
&&\int_{\textbf{X}}e^{\tilde{l}(x,\tilde{a}^\ast,y)}W^\ast(y)\tilde{Q}(dy|x,\tilde{a}^\ast)=\inf_{\tilde{a}\in\tilde{\textbf{A}}}\left\{\int_{\textbf{X}}e^{\tilde{l}(x,\tilde{a},y)}W^\ast(y)\tilde{Q}(dy|x,\tilde{a})\right\}\nonumber\\
&=&\int_{\textbf{X}}e^{\tilde{l}(x,\psi^\ast(x),y)}W^\ast(y)\tilde{Q}(dy|x,\psi^\ast(x))I\{\tilde{a}^\ast\in\textbf{A}^I\}+
\int_{\textbf{X}}e^{\tilde{l}(x,f^\ast(x),y)}W^\ast(y)\tilde{Q}(dy|x,f^\ast(x))I\{\tilde{a}^\ast\in\textbf{A}^G\}.\nonumber\\
\end{eqnarray}
Such a pair $(\psi^\ast,f^\ast)$ of measurable selectors exists.
\end{proposition}

We introduce the notation to be used in the next statement. Define for each $[1,\infty]$-valued universally measurable function $g$ on $\textbf{X}$
\begin{eqnarray}\label{JapanRevisionAAEqn02}
{\textbf{X}^G(g):=\left\{x\in\textbf{X}:\infty>g(x),~0=\inf_{a\in \textbf{A}^G}\left\{\int_{\textbf{X}}g(y)\tilde{q}(dy|x,a)-(q_x(a)-c^G(x,a))g(x)\right\}\right\}},
\end{eqnarray}
and denote by $\textbf{X}^I(g)$ the collection of $x\in\textbf{X}$ at which, $g(x)=\inf_{b\in\textbf{A}^I}\left\{\int_{\textbf{X}}g(y)e^{c^I(x,b,y)}Q(dy|x,b)\right\}.$ {Proposition \ref{GGZyExponentialProposition02} in the Appendix asserts that $W^\ast$ is universally measurable so that the integrals $\int_{\textbf{X}}W^\ast(y)\tilde{q}(dy|x,a)$ and $\int_{\textbf{X}}W^\ast(y)e^{c^I(x,b,y)}Q(dy|x,b)$ are defined.}

{
\begin{theorem}\label{JapanTheorem05}
Suppose Conditions \ref{JapanCon02} and \ref{JapanCon01} are satisfied. Then the following assertions hold.
\begin{itemize}
\item[(a)] The value function $W^\ast$ of problem (\ref{ZyExponentialProblem2}) for the tilde DTMDP model coincides with ${\cal V}^\ast$.
\item[(b)] $\textbf{X}\setminus \textbf{X}^I(W^\ast)\subseteq \textbf{X}^G(W^\ast).$
\item[(c)] There is a deterministic stationary optimal policy for the gradual-impulse control problem (\ref{JapanProblem01}), which can be obtained as follows. For each pair $(\psi^\ast,f^\ast)$ of measurable mappings satisfying (\ref{JapanConserving01}) (and there exists such a pair by Proposition \ref{JapanRevisionAAProp01}), the following deterministic stationary policy $(\varphi,\psi,F)$ is optimal, where
\begin{eqnarray*}
\psi(x)=\psi^\ast(x),~F(x)_t(da)\equiv \delta_{f^\ast(x)}(da)
\end{eqnarray*}
for all $x\in\textbf{X}$, and $\varphi(x)=\infty$ (respectively, $\varphi(x)=0$) for all $x\in \textbf{X}\setminus \textbf{X}^I(W^\ast)$ (respectively $x\in\textbf{X}^I(W^\ast)$).
\end{itemize}
\end{theorem}}
The proofs of this and the other statements in this section are postponed to Section \ref{JapanSecProof}.

According to Theorem \ref{JapanTheorem05}, roughly speaking, if the current state is in $\textbf{X}^G(W^\ast)$, then it is optimal not to apply impulse until the next natural jump; and if the current state is in $\textbf{X}^I(W^\ast)$, then it is optimal to apply immediately an impulse.

{
According to Theorem \ref{JapanTheorem05}, (\ref{JapanBellman05}) is the optimality equation for the gradual-impulse control problem (\ref{JapanProblem01}). It can be written out in an equivalent form that does not involve the function $w$, which might be more convenient sometimes.
\begin{corollary}\label{JapanRevisionCorollary0001}
Suppose Conditions \ref{JapanCon02} and \ref{JapanCon01} are satisfied. Then the following assertions hold.
\begin{itemize}
\item[(a)] ${\cal V}^\ast$ is the minimal $[1,\infty]$-valued lower semicontinuous function on $\textbf{X}$ satisfying
\begin{eqnarray}\label{JapanREqnA}
&&\inf_{a\in\textbf{A}^G}\left\{\int_{\textbf{X}} {\cal V}^\ast(y) \tilde{q}(dy|x,a)-(q_x(a)-c^G(x,a)){\cal V}^\ast(x)  \right\}\ge 0,\\
&&~\forall~x\in\textbf{X}^\ast({\cal V}^\ast):=\{x\in\textbf{X}:~{\cal V}^\ast(x)<\infty\}\nonumber
\end{eqnarray}
and
\begin{eqnarray}\label{JapanREqnB}
{\cal V}^\ast(x)\le \inf_{b\in\textbf{A}^I}\left\{ \int_{\textbf{X}} e^{c^I(x,b,y)}{\cal V}^\ast(y)Q(dy|x,b) \right\},~x\in\textbf{X},
\end{eqnarray}
whereas at each $x\in\textbf{X}$,  the inequality in either (\ref{JapanREqnA}) or (\ref{JapanREqnB}) holds with equality.
\item[(b)] A pair  $(\psi^\ast,f^\ast)$ of measurable mappings satisfies (\ref{JapanConserving01}) if and only if
\begin{eqnarray*}
&&\inf_{a\in \textbf{A}^G}\left\{\int_{\textbf{X}}{\cal V}^\ast(y)\tilde{q}(dy|x,a)-(q_x(a)-c^G(x,a)){\cal V}^\ast(x)\right\}\\
&=&\int_{\textbf{X}}{\cal V}^\ast(y)\tilde{q}(dy|x,f^\ast(x))-(q_x(f^\ast(x))-c^G(x,f^\ast(x))){\cal V}^\ast(x)
\end{eqnarray*}
for each $x\in\textbf{X}^G({\cal V}^\ast)$, and
\begin{eqnarray*}
\inf_{b\in\textbf{A}^I}\left\{\int_{\textbf{X}}e^{c^I(x,b,y)}{\cal V}^\ast(y)Q(dy|x,b)\right\}=\int_{\textbf{X}}{\cal V}^\ast(y)e^{c^I(x,\psi^\ast(x),y)}Q(dy|x,\psi^\ast(x)),~\forall~x\in\textbf{X}.
\end{eqnarray*}
(According to Theorem \ref{JapanTheorem05}, $(\psi^\ast,f^\ast)$ gives rise to a deterministic stationary optimal policy for the gradual-impulse control problem (\ref{JapanProblem01}).)
%\item[(d)] For each pair of measurable mappings $(\psi^\ast,f^\ast)$ that satisfy the previous two relations,
\end{itemize}
\end{corollary}
}

{Under the conditions of the previous statement, in the first glance, given ${\cal V}^\ast$ being an $[1,\infty]$-valued lower semicontinuous function on $\textbf{X}$, it may be not immediately clear why the claimed measurable selector $f^\ast$ exists because in \begin{eqnarray*}
\int_{\textbf{X}}{\cal V}^\ast(y)\tilde{q}(dy|x,a)-(q_x(a)-c^G(x,a)){\cal V}^\ast(x)=\left(\int_{\textbf{X}}{\cal V}^\ast(y)\tilde{q}(dy|x,a)+c^G(x,a){\cal V}^\ast(x)\right)-(q_x(a){\cal V}^\ast(x))
\end{eqnarray*}
the expressions in the two brackets are both lower semicontinuous in $(x,a)\in\textbf{X}\times\textbf{A}^G$, and the difference between two lower semicontinuous functions may be not lower semicontinuous. This and
Lemma \ref{JapanLemma05} are the motivation of considering the tilde DTMDP model.} %We shall show the lower semicontinuity of $V^\ast$ and the existence of the required selectors by reducing the original model to the following simpler DTMDP, which will be referred to as the tilde DTMDP model.

To end this section, we present a simple example to demonstrate a situation, where it is natural and necessary to allow multiple impulses at a single time moment.
\begin{example}\label{JapanRevisionExample005}
Let us revisit Example \ref{JapanRevisionExample001}. The model has a state space $\{1,2\}$, where $1$ stands for the rat being present in the kitchen, and $2$ indicates the rat either dead or outside the house. The space of gradual controls is a singleton and will not be indicated explicitly, and the space of impulses is $\textbf{A}^I=\{0,1\}$, with $1$ or $0$ standing for shooting or not. So the inequalities (\ref{JapanREqnA}) and (\ref{JapanREqnB}) for the value function ${\cal V}^\ast$ read:
\begin{eqnarray*}
&&{\cal V}^\ast(2)=1;~\mu {\cal V}^\ast(2) -(\mu-l){\cal V}^\ast(1)\ge 0;~{\cal V}^\ast(1)\le \min\{e^C p {\cal V}^\ast(2)+e^C (1-p){\cal V}^\ast(1),{\cal V}^\ast(1)\}.
\end{eqnarray*}
Suppose $1-e^C(1-p)>0.$ By Theorem \ref{JapanTheorem05} and Corollary \ref{JapanRevisionCorollary0001}, if $\frac{e^C p}{1-e^C(1-p)}>\frac{\mu}{\mu-l}>0$, then ${\cal V}^\ast(1)=\frac{\mu}{\mu-l}$, and the optimal deterministic stationary policy is to never shoot at the rat; otherwise, ${\cal V}^\ast(1)=\frac{e^C p}{1-e^C(1-p)}=E[e^{CZ}]$ with $Z$ following the geometric distribution with success probability $p$, and the optimal deterministic stationary policy is to keep shooting as soon as the rat is in kitchen until the rat was hit.
\end{example}

The proofs of the statements in this section are based on the investigation of an optimal control problem for another DTMDP model, which will be referred to as the hat DTMDP model and introduced in Section \ref{JapanSecDTMDPhat}. For this moment, we point out that the hat DTMDP model is quite different from the tilde DTMDP model: it is with a more complicated action space, and is not necessarily  semicontinuous under Conditions \ref{JapanCon02} and \ref{JapanCon01}, see Examples \ref{JapanRevisionExample002} and \ref{JapanExample01}.

\section{The hat DTMDP model}\label{JapanSecDTMDPhat}

\begin{figure}
  \includegraphics[width=1.0\textwidth, angle = 0]{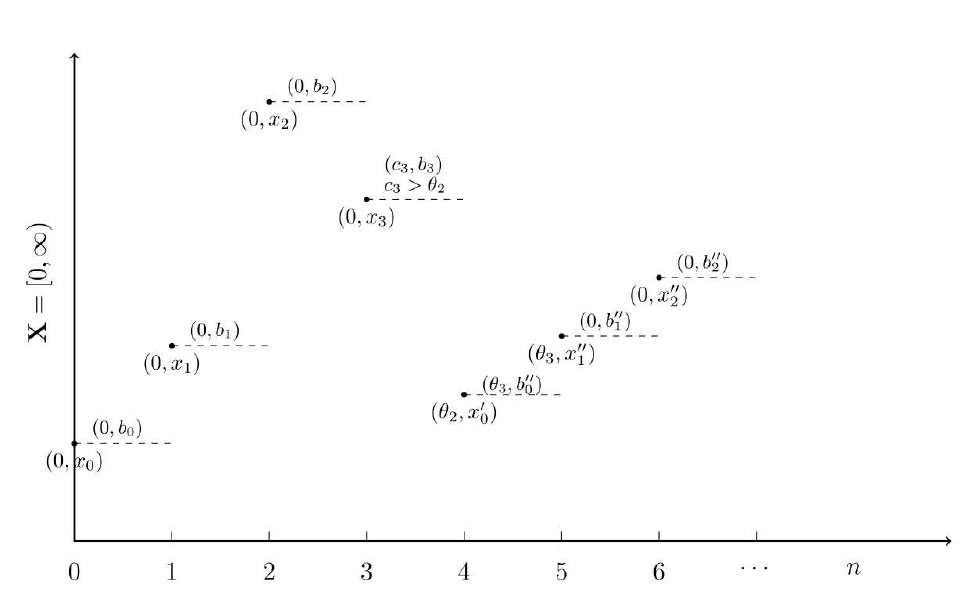}
  \caption{\small {The realization of the state process in the hat DTMDP model corresponding to sample path in the gradual-impulse control in Figure \ref{Japanfig:01}. The time index is discrete from $\{0,1,\dots\}$. The realizations of the components $\{(C_n,B_n)\}_{n=0}^\infty$ in the action process $\{\hat{A}_n\}_{n=0}^\infty$ are indicated above the dashed lines between consecutive states. For example, $(0,b_0)$ next to the state $(0,x_0)$ indicates that the decision maker applies an impulse $b_0$ immediately, which results in the next state $(0,x_1).$ All the components $x_0,x_1,\dots,$ $x_0'$, $x_1''$, $x_2''$ and $b_1,b_2,b_0'',b_1'',b_2''$ are the same as in Figure \ref{Japanfig:01}. The only exception is $(c_3,b_3),$ which does not appear in Figure \ref{Japanfig:01}. Nevertheless, $c_3>\theta_2$, because in Figure \ref{Japanfig:01}, the first jump in the marked point process therein at the time moment $\theta_1+\theta_2=\theta_2$ is triggered by a natural jump.}}
  \label{Japanfig:02}
\end{figure}

In this section, we describe a DTMDP problem, which will serve the investigations of the gradual-impulse control problem.  To distinguish it from the intervention DTMDP model, we shall refer to it as the hat DTMDP model. {The system primitives of the DTMDP model are defined in terms of those of the gradual-impulse control problem.} We will reveal, in greater detail, the connections relevant to this paper between the hat DTMDP problem and the gradual-impulse control problem at the end of this section. {For a first impression, roughly speaking, the state process of the hat DTMDP model comes from the system dynamics of the gradual impulse control problem in the following way. The state has two coordinates. Along the (discrete-time) state process of the hat DTMDP model, the second coordinates record the system states of the graduate-impulse control problem immediately after a natural jump (of the marked point process $\{(T_n,Y_n)\}_{n=1}^\infty$) or an ``actual'' impulse (thus the state immediately after the psuedo impulse $\Delta$ will not be recorded). The first coordinates record the time in the gradual-impulse control problem elapsed between two consecutive states as recorded in the second coordinates. For the sake of illustration, the realization of the state process in the hat DTMDP model corresponding to the sample path in Figure \ref{Japanfig:01} of the gradual-impulse control problem  is displayed in Figure \ref{Japanfig:02}.}

The hat DTMDP is with a more complicated action space as compared with the original gradual-impulse control problem. {To describe the action space of the hat DTMDP model, let us recall some known and general facts and introduce some further notations as follows.}

 {Let ${\cal R}$ be the collection of ${\cal P}(\textbf{A}^G)$-valued measurable mappings on $[0,\infty)$ with any two elements therein being identified the same if they differ only on a null set with respect to the Lebesgue measure. Recall that ${\cal P}(\textbf{A}^G)$ stands for the space of probability measures on $(\textbf{A}^G, {\cal B}(\textbf{A}^G))$. We endow ${\cal P}(\textbf{A}^G)$ with its weak topology (generated by bounded continuous functions on $\textbf{A}^G$) and the Borel $\sigma$-algebra, so that ${\cal P}(\textbf{A}^G)$ is a Borel space, see Chapter 7 of \cite{Bertsekas:1978}. It is known, see Lemma 1 of \cite{Yushkevich:1980}, that the space ${\cal R}$, endowed with the smallest $\sigma$-algebra with respect to which the mapping $\rho=(\rho_t(da))\in{\cal R}\rightarrow \int_0^\infty e^{-t}g(t,\rho_t)dt$ is measurable for each bounded measurable function $g$ on $(0,\infty)\times {\cal P}(\textbf{A}^G)$, is a Borel space.
Then, according to Section 43 of \cite{Davis:1993}, the space ${\cal R}$ is a compact metrizable space, endowed with the Young topology, which is the coarsest topology with respect to which, the mapping
\begin{eqnarray*}
\rho=(\rho_t(da))\in {\cal R}\rightarrow \int_0^\infty \int_{\textbf{A}^G} g(t,a)\rho_t(da)dt
\end{eqnarray*}
is continuous for each function $g$ on $(0,\infty)\times \textbf{A}^G$ satisfying that (a) for each $t\in(0,\infty),$ $g(t,\cdot)$ is continuous on $\textbf{A}^G$; (b) for each $a\in\textbf{A}^G$, $g(\cdot,a)$ is measurable on $(0,\infty);$ and (c) $\int_0^\infty \sup_{a\in\textbf{A}^G}|g(t,a)|dt<\infty.$ Such a function $g$ satisfying these requirements is called a strongly integrable Caratheodory function.}

{
Below we shall use, without special reference, the following notation. If $\mu$ is a measure on a Borel space $(X,{\cal{B}}(X))$, then the notation $f(\mu):=\int_X f(x)\mu(dx)$ is in use for each measurable function $f$ on $(X,{\cal{B}}(X))$, provided that the integral is well defined.}

%In the next section, we will reduce the gradual-impulse control problem and the hat DTMDP problem to yet another simpler DTMDP problem.

\subsection{Primitives of the hat DTMDP model}
The state space of the hat DTMDP model is $\hat{\textbf{X}}:=\{(\infty,x_\infty)\}\bigcup [0,\infty)\times\textbf{X}$, where $(\infty,x_\infty)$ is an isolated point, and the action space of the DTMDP is $\hat{\textbf{A}}:=[0,\infty]\times \textbf{A}^I\times {\cal R}$. Endowed with the product topology, where $[0,\infty]$ is compact in the standard topology of the extended real-line, $\hat{\textbf{A}}$ is also a compact Borel space. Here, $\textbf{X}$, $\textbf{A}^I$ and $\textbf{A}^G$ are the state, impulse and gradual action spaces in the gradual-impulse control problem.

The transition probability $p$ is defined as follows, where the notation introduced above this subsection is in use, e.g., $q_x(\rho_t):=\int_{\textbf{A}^G} q_x(a)\rho_t(da)$ and $c^G(x,\rho_t):=\int_{\textbf{A}^G}c^G(x,a)\rho_t(da)$.
For each bounded measurable function $g$ on $\hat{\textbf{X}}$ and action $\hat{a}=(c,b,\rho)\in\hat{\textbf{A}}$,
\begin{eqnarray*}
&&\int_{\hat{\textbf{X}}}g(t,y)p(dt\times dy|(\theta,x),\hat{a})\\
&:=&I\{c=\infty\}\left\{ g(\infty,x_\infty)e^{-\int_0^\infty q_x(\rho_s)ds}+\int_{0}^\infty\int_{\textbf{X}}g(t,y)\tilde{q}(dy|x,\rho_t)e^{-\int_0^t q_x(\rho_s)ds} dt  \right\}\\
&&+I\{c<\infty\}\left\{\int_0^c \int_{\textbf{X}} g(t,y)\tilde{q}(dy|x,\rho_t)e^{-\int_0^t q_x(\rho_s)ds}dt+e^{-\int_0^c q_x(\rho_s)ds}\int_{\textbf{X}}g(c,y)Q(dy|x,b)  \right\}\\
&=&\int_0^c \int_{\textbf{X}} g(t,y)\tilde{q}(dy|x,\rho_t)e^{-\int_0^t q_x(\rho_s)ds}dt+I\{c=\infty\} g(\infty,x_\infty)e^{-\int_0^\infty q_x(\rho_s)ds}\\
&&+I\{c<\infty\}e^{-\int_0^c q_x(\rho_s)ds}\int_{\textbf{X}}g(c,y)Q(dy|x,b)
\end{eqnarray*}
for each state
  $(\theta,x)\in[0,\infty)\times \textbf{X}$; and
  \begin{eqnarray*}
&&\int_{\hat{\textbf{X}}}g(t,y)p(dt\times dy|(\infty,x_\infty),\hat{a}):=g(\infty,x_\infty).
  \end{eqnarray*}

It is known, see e.g., \cite{Costa:2013,Forwick:2004}, that for each bounded measurable function $g$ on $\hat{\textbf{X}}$, the above expressions are indeed measurable on $\hat{\textbf{X}}\times \hat{\textbf{A}}$, and the same also concerns the cost function $l$ on $\hat{\textbf{X}}\times\hat{\textbf{A}}\times\hat{\textbf{X}}$ defined as follows:
\begin{eqnarray*}
l((\theta,x),\hat{a},(t,y)):=I\{(\theta,x)\in [0,\infty)\times\textbf{X}\}\left\{\int_0^t c^G(x,\rho_s)ds+I\{t=c\}c^I(x,b,y)\right\}
\end{eqnarray*}
for each  $(\theta,x),\hat{a},(t,y))\in\hat{\textbf{X}}\times\hat{\textbf{A}}\times\hat{\textbf{X}}$, accepting that $c^I(x,b,x_\infty)\equiv 0.$ Recall that the generic notation $\hat{a}=(c,b,\rho)\in\hat{\textbf{A}}$ of an action in this hat DTMDP model has been in use. The pair $(c,b)$ is the pair of the planned time until the next impulse and the next planned impulse, and $\rho$ is (the rule of) the relaxed control to be used during the next sojourn time. {The realization of the components $\{(C_n,B_n)\}_{n=0}^\infty$ of the action process in the hat DTMDP model corresponding to the sample path in Figure \ref{Japanfig:01} of the gradual-impulse control problem is displayed in Figure \ref{Japanfig:02}.}

{
For the convenience in future reference, we make the following definition.
\begin{definition}
The hat DTMDP model is the following four-tuple $\{\hat{\textbf{X}}, \hat{\textbf{A}}, p, l\}$, all defined above in terms of the primitives of the gradual-impulse control problem.
\end{definition}
}

Note that Condition \ref{JapanCon01} does not imply that the hat DTMDP model is semicontinuous, which is defined in the appendix. In fact, the transition probability $p$, in general, does not satisfy the weak continuity condition, even under Condition \ref{JapanCon01}. This is demonstrated by the next two examples.
\begin{example}\label{JapanRevisionExample002}
Suppose $q_x(a)\equiv 0$, and $\textbf{A}^G$ and $\textbf{A}^I$ are both singletons. Consider $\hat{a}_n=(c_n,b,\rho)$, where $c_n\rightarrow \infty$ and $c_n\in[0,\infty)$ for each $n\ge 1$; and the bounded continuous function on $\hat{\textbf{X}}$: $g(t,x)\equiv 1$ for each $(t,x)\in[0,\infty)\times \textbf{X}$, and $g(\infty,x_\infty)=0$. Then $\int_{\hat{\textbf{X}}}g(t,y)p(dt\times dy|(\theta,x),\hat{a}_n)=\int_{\textbf{X}}g(c_n,y)Q(dy|x,b)=1$ for each $n\ge 1,$ whereas $\int_{\hat{\textbf{X}}}g(t,y)p(dt\times dy|(\theta,x),(\infty,b,\rho))=g(\infty,x_\infty)=0\ne 1$.
\end{example}

\begin{example}\label{JapanExample01}
Consider $\textbf{A}^G=[0,1],$ $\textbf{A}^I$ an arbitrary compact Borel space, $\textbf{X}$ a finite set (endowed with discrete topology),   $q_x(a)=a$ for each $x\in\textbf{X}$. Then consider $x^{(n)}\equiv x\in \textbf{X}$, $b^{(n)}\equiv b$, $c^{(n)}\equiv c=\infty,$ and for each $t\ge 0,$ $\rho^{(n)}_t(da)=\delta_{\frac{1}{n}}(da)$, and $\rho_t(da)=\delta_0(da).$
Then for each strongly integrable Caratheodory function $g(t,a)$,
\begin{eqnarray*}
\int_0^\infty g(t,\rho^{(n)}_t)dt-\int_0^\infty g(t,\rho^{(0)}_t)dt=\int_0^\infty (g(t,\frac{1}{n})-g(t,0))dt\rightarrow 0
\end{eqnarray*}
as $n\rightarrow \infty,$ by using the dominated convergence theorem. Thus, $\rho^{(n)}\rightarrow \rho$ as $n\rightarrow \infty.$ Let $\hat{a}_n=(c^{(n)},b^{(n)},\rho^{(n)})$ and $\hat{a}=(c,b,\rho)$. It follows that $ ((\theta,x),\hat{a}_n)\rightarrow ((\theta,x),\hat{a})$ as $n\rightarrow \infty.$
Now consider the bounded continuous function on $\hat{\textbf{X}}$ given by $g(\infty,x_\infty)= 1$ and $g(t,x)\equiv 0$ on $[0,\infty)\times {\textbf{X}}$. (Recall that $(\infty,x_\infty)$ is an isolated point in $\hat{\textbf{X}}$.) Then we see
\begin{eqnarray*}
&&\lim_{n\rightarrow \infty}\int_{\hat{\textbf{X}}}g(t,y)p(dt\times dy|(\theta,x),\hat{a}_n)=\lim_{n\rightarrow \infty}e^{-\int_0^\infty q_x(\rho^{(n)}_s)ds}=\lim_{n\rightarrow \infty}e^{-\int_0^\infty \frac{1}{n}ds}=0\\
&\ne& 1=e^{-\int_0^\infty 0 ds}=\int_{\hat{\textbf{X}}}g(t,y)p(dt\times dy|(\theta,x),\hat{a}).
\end{eqnarray*}
\end{example}

\begin{remark}
Example \ref{JapanExample01} implies that the assertion of Lemma 5.12 in \cite{Zhang:2017} (stated without proof) is inaccurate without further conditions (such as $q_x(a)>\delta>0$ for some $\delta>0$). Similarly, Lemma 4.1(b) in \cite{XinGuo:2018} is correct if $q_x(a)>\delta>0$ for some $\delta>0$.  However, the optimality results in \cite{Zhang:2017} all survive without assuming extra conditions, as a particular consequence of the arguments presented below in the present paper. %The same remark also applies to the fixing of the inaccurate statement of Lemma 4.1(b) in \cite{XinGuo:2018}. %We mention that the inconvenient term is the first expression in $F(x,\hat{a})$, and it does not appear in the risk neutral case or if $q_x(a)>\delta>0$ for each $(x,a)\in\textbf{X}\times\textbf{A}^G$, c.f. Condition \ref{JapanAssumption01} and the observation after it.
\end{remark}

We use the notation $\hat{h}_n=((\theta_0,x_0),(c_0,b_0,\rho_0),(\theta_1,x_1), (c_1,b_1,\rho_1),(\theta_2,x_2),\dots,(\theta_n,x_n))$ for the $n$-history in the hat DTMDP model.

%{In a trajectory of system dynamics of the gradual-impulse control problem, when a sequence of impulses are applied in turn at a single time moment, they result in a sequence of states displayed vertically at that single time moment, see Figure. Roughly speaking, the state process in the hat DTMDP model is obtained when we view such ``vertical'' sequence horizontally. The second coordinate of a state in the hat DTMDP model records the system state in the gradual-impulse control problem, whereas the first coordinate records the time elapsed between the current state (given by the second coordinate) and the previous system states.
%}

The concerned optimal control problem for the hat DTMDP model reads:
\begin{eqnarray}\label{JapanProblem02}
\mbox{Minimize over $\sigma$:}~\mathbb{E}_{\hat{x}}^\sigma\left[e^{\sum_{n=0}^\infty l(\hat{X}_n,\hat{A}_n,\hat{X}_{n+1})}\right]=:V((\theta,x),\sigma)
\end{eqnarray}
where $\{\hat{X}_n\}_{n=0}^\infty$ and $\{\hat{A}_n\}_{n=0}^\infty$ are the state and action processes, and the minimization problem is over all strategies $\sigma$ in the hat DTMDP model. (See the appendix for the basic notations in a DTMDP.)
We denote by $V^\ast$ the value function of this optimal control problem, i.e.,
\begin{eqnarray*}
V^\ast(\theta,x):=\inf_{\sigma}\mathbb{E}_{\hat{x}}^\sigma\left[e^{\sum_{n=0}^\infty l(\hat{X}_n,\hat{A}_n,\hat{X}_{n+1})}\right]
\end{eqnarray*}
for each $\hat{x}=(\theta,x)\in\hat{\textbf{X}}$, where the infimum is over all strategies. Clearly, $V^\ast(\infty,x_\infty)=1$. It will be seen in Lemma \ref{JapanNewLem01} that $V^\ast$ depends on $(\theta,x)$ only through $x$, and a strategy $\sigma$ is optimal if
\begin{eqnarray*}
V((0,x),\sigma)=V^\ast(x)
\end{eqnarray*}
for each $x\in\textbf{X}$.  Below, when the context is clear, we often consider the restriction of $V^\ast$ on $\textbf{X}$ but still use the same notation. The definition of an optimal strategy and other relevant notions of DTMDP are collected in the appendix.

%Let us say a few words regarding some notations and conventions to be used below for brevity. Let $\hat{h}_n=((\theta_0,x_0),(c_0,b_0,\rho_0),(\theta_1,x_1), (c_1,b_1,\rho_1),(\theta_2,x_2),\dots,(\theta_n,x_n))$ be the $n$-history in the hat DTMDP model.
Consider a strategy $\sigma=\{\sigma_n\}_{n=0}^\infty$ in the hat DTMDP model, where for each $n\ge 0,$ $\sigma_n(d\hat{a}|\hat{h}_n)$ is a stochastic kernel on $\hat{\textbf{A}}$ given $\hat{h}_n$, which specifies the conditional distribution of the next action $(c,b,\rho)$ given $\hat{h}_n$.

In general, a strategy in the hat DTMDP model can make use of past decision rules of relaxed controls, and the selection of the next relaxed control, and that of the next planned impulse time and impulse do not have to be (conditionally) independent. Therefore, a general strategy in the hat DTMDP model does not immediately correspond to a policy in the continuous-time gradual-impulse control problem described in the previous section. To relate the continuous-time gradual-impulse control problem (\ref{JapanProblem01}) and the hat DTMDP problem (\ref{JapanProblem02}), see Proposition \ref{JapanProposition01} below, we introduce the following class of strategies in the hat DTMDP model.

\begin{definition}
A strategy $\sigma$ in the hat DTMDP model is called typical if under it, given $\hat{h}_n$, the selection of the next action $(c,b)$ and $\rho$ are conditionally independent,  and {moreover}, the selection of $\rho$ is deterministic, i.e., \begin{eqnarray*}
\sigma_n(dc\times db\times d\rho|\hat{h}_n)=\sigma_n'(dc\times db|\hat{h}_n)\delta_{F^n(\hat{h}_n)}(d\rho),
\end{eqnarray*}
where $F^n(\hat{h}_n)$ is measurable in its argument and takes values in ${\cal R}$, and  $\sigma_n'(dc\times db|\hat{h}_n)$ is a stochastic kernel on $[0,\infty]\times \textbf{A}^I$ given $\hat{h}_n$.
\end{definition}

One can always write  $\sigma_n'(dc\times db|\hat{h}_n)=\varphi_n(dc|\hat{h}_n)\psi_n(db|\hat{h}_n,c)$ for some stochastic kernels $\varphi_n$ and $\psi_n$. Intuitively, $\varphi_n$ defines the (conditional) distribution of the planned time duration till the next impulse, and $\psi_n(db|\hat{h}_n,c)$ specifies the distribution of the next impulsive action given  the history $\hat{h}_n$ and the next impulse moment $c$, provided that it takes place before the next natural jump.  Therefore, we identify a typical strategy $\sigma=\{\sigma_n\}$ as $\{(\varphi_n,\psi_n,F^n)\}_{n=0}^\infty.$

For further notational brevity, when the stochastic kernels $\varphi_n$ are identified with underlying measurable mappings, we will use $\varphi_n$ for the measurable mappings, and write $\varphi_n(\hat{h}_n)$ instead of $\varphi_n(da|\hat{h}_n)$. The same applies to other stochastic kernels such as $\psi_n$. The context will exclude any potential confusion.

Finally, in general, we often do not indicate the arguments that do not affect the values of the concerned mappings.  For example, if $\varphi_n(\hat{h}_n)$ depends on $\hat{h}_n$ only through $x_n$, then we write $\varphi_n(da|\hat{h}_n)$ as $\varphi_n(da|x_n)$.

\subsection{Connection between the gradual-impulse control problem and the hat DTMDP problem}\label{JapanRevisionSubsection001}
Each policy $u$ as introduced in Definition \ref{JapanDef01} induces a (typical) strategy $\{(\varphi_n,\psi_n,F^n)\}_{n=0}^\infty$ in the hat DTMDP model as follows, where we only need consider $x_n\in\textbf{X}$, as the definition of the strategies at $x_n=x_\infty$ is immaterial, and can be arbitrary. For each $m\ge 1$, and $h_m\in\textbf{H}_m$, there exists a strategy $\pi^{\Gamma^1_m,h_m}=\{\pi^{\Gamma^1_m,h_m}_n\}_{n=0}^\infty$ in the intervention DTMDP model such that $\Gamma^1_m(dy|h_m)=\beta^{\pi^{\Gamma^1_m,h_m}}(dy|\bar{x}(y_m)).$ Similarly, for each $x\in\textbf{X},t>0$,
there exists a strategy $\pi^{\Gamma^0_m,h_m,t,x}=\{\pi^{\Gamma^0_m,h_m,t,x}_n\}_{n=0}^\infty$ in the intervention DTMDP model such that $\Gamma^0_m(dy|h_m,t,x)=\beta^{\pi^{\Gamma^0_m,h_m,t,x}}(dy|x).$ Finally,  there is a strategy
$\pi^{u_0}=\{(\pi^{u_0}_n)\}_{n=0}^\infty$
in the intervention DTMDP model satisfying \begin{eqnarray}\label{JapanRevision02}
u_0(dy|x)=\beta^{\pi^{u_0}}(dy|x)
\end{eqnarray}
for each $x\in\textbf{X}$.

Consider the case of $n=0$. Then we define
\begin{eqnarray*}
\varphi_0(\{0\}|\theta,x)&:=& 1-\pi_0^{u_0}(\{\Delta\}|x);\\
\varphi_0(dc|\theta,x)&:=&\pi_0^{u_0}(\{\Delta\}|x)\Phi_1(dc|(x,\Delta,\Delta,\dots),0,(x,\Delta,\Delta,\dots)) ~\mbox{on~} (0,\infty];\\
\psi_0(db|\theta,x,c)&:=&\frac{\pi^{u_0}_0(db|x)}{1-\pi^{u_0}_0(\{\Delta\}|x)}I\{c=0\}+I\{c>0\} \frac{\pi^{\Gamma^1_1,((x,\Delta,\dots),0,(x,\Delta,\dots))}_0(db|x)}{1-\pi^{\Gamma^1_1,((x,\Delta,\dots),0,(x,\Delta,\dots))}_0(\{\Delta\}|x)}\\
&=&\frac{\pi^{u_0}_0(db|x)}{1-\pi^{u_0}_0(\{\Delta\}|x)}I\{c=0\}+I\{c>0\}  \pi^{\Gamma^1_1,((x,\Delta,\dots),0,(x,\Delta,\dots))}_0(db|x);\\
F^0(\theta,x)_t(da)&:=&\Pi_1(da|(x,\Delta,\Delta,\dots),0,(x,\Delta,\Delta,\dots),t),
\end{eqnarray*}
where the second equality in the definition of $\psi_0(db|\theta,x,c)$ holds because $\pi^{\Gamma^1_1,((x,\Delta,\dots),0,(x,\Delta,\dots))}_0(\{\Delta\}|x)=0$, which follows from the requirement that $\Gamma^1_n(\cdot|h_n)\in {\cal P}^{\textbf{Y}^\ast}(\bar{x}(y_n))$ for all $n\ge 1$ in Definition \ref{JapanDef01}. Also concerning the definition of $\psi_0(db|\theta,x,c)$, note that if the denominator in $ {1-\pi^{u_0}_0(\{\Delta\}|x)}=0,$ we put $\frac{\pi^{u_0}_0(db|x)}{1-\pi^{u_0}_0(\{\Delta\}|x)}$ as an arbitrary stochastic kernel. The reason is that in the expression $\frac{\pi^{u_0}_0(db|x)}{1-\pi^{u_0}_0(\{\Delta\}|x)}I\{c=0\}$, equality ${1-\pi^{u_0}_0(\{\Delta\}|x)}=0$ would indicate that the probability of selecting an instantaneous impulse is zero, and so $I\{c=0\}=0$ almost surely. The same explanation applies to the definitions of $\psi_n(db|\hat{h}_n,c)$ below, and will not be repeated there. Note that the right hand side does not depend on $\theta\in[0,\infty)$, because the initial time moment is always fixed to be $\theta=0$.

 {The intuition behind the above definition of $(\varphi_0,\psi_0,F^0)$ is as follows. Recall that, if the initial system state is $x\in\textbf{X}$, then the intervention $y_1\in\textbf{Y}$ at the initial time in the gradual-impulse control problem is a realization from the distribution $u_0(\cdot|x)=\beta^{\pi^{u_0}}(\cdot|x),$ which is the strategic measure of some strategy $\pi^{u_0}=\{\pi_n^{u_0}\}_{n=0}^\infty$ in the intervention DTMDP model, see (\ref{JapanRevision02}). Then $\pi_0^{u_0}(\{\Delta\}|x)$ is the probability that no impulse is applied at the initial time $0$ (given the initial system state $x$) in the gradual-impulse control problem. Consequently, $1-\pi_0^{u_0}(\{\Delta\}|x)$ is the probability to apply an impulse immediately, i.e., to wait time $0$ until the next impulse, and thus $\varphi_0(\{0\}|\theta,x)$. This quantity does not depend on $\theta$, because the initial time is always $0.$ Then for a measurable subset $\Gamma_1\subseteq(0,\infty]$,
\begin{eqnarray*}
&&\pi_0^{u_0}(\{\Delta\}|x)\Phi_1(\Gamma_1|(x,\Delta,\Delta,\dots),0,(x,\Delta,\Delta,\dots))\\
&=&\mbox{Probability (no impulse at initial time $0$ given initial system state $x$)}\\
&&\times \mbox{Probability (time to wait until next impulse is in $\Gamma_1$} \\
&&\mbox{given no impulse is immediately applied at the initial time with the initial state $x$)},
\end{eqnarray*}
which is equal to
\begin{eqnarray*}
&&\mbox{Probability (No immediate impulse, and the time duration until the next planned impulse is in $\Gamma$)}\\
&=&\mbox{Probability (the time duration until the next planned impulse is in $\Gamma$)},
\end{eqnarray*}
and thus $\varphi_0(\Gamma|\theta,x)$, where the equality follows because $\Gamma\subseteq(0,\infty].$ (Recall that a planned impulse takes place if no natural jump occurs during the time duration to wait for it.) Finally, as for $\psi_0(db|\theta,x,c)$, if $c=0$, and $\Gamma_2\in{\cal B}(\textbf{A}^I)$, then
\begin{eqnarray*}
&&\frac{\pi^{u_0}_0(\Gamma_2|x)}{1-\pi^{u_0}_0(\{\Delta\}|x)}=\frac{\mbox{Probability (an immediate impulse from $\Gamma_2$ is applied)}}{\mbox{Probability(an immediate impulse is applied)}}\\
&=&\mbox{Probability (an impulse is applied immediately from $\Gamma_2$}\\
&&\mbox{given that an impulse is applied after time duration $0$)},
\end{eqnarray*}
which is thus $\psi_0(\Gamma_2|\theta,x,0).$ One can understand $\psi_0(db|\theta,x,c)$ when $c>0$ in the same manner. The very similar intuition guides the definition of $(\varphi_n,\psi_n,F^n)$ below.}

Now consider $n\ge 1.$ Let $\hat{h}_n=((\theta_0,x_0),(c_0,b_0,\rho_0),(\theta_1,x_1), (c_1,b_1,\rho_1),(\theta_2,x_2),\dots,(\theta_n,x_n))$ be the $n$-history in the hat DTMDP model.
If $\{1\le i\le n: ~\theta_i>0\}= \emptyset,$ then we define
\begin{eqnarray*}
&&\varphi_n(\{0\}|\hat{h}_n):=1-\pi^{u_0}_n(\{\Delta\}|x_0,b_0,\dots,b_{n-1},x_n),\\
&&\varphi_n(dc|\hat{h}_n):=\pi^{u_0}_n(\{\Delta\}|x_0,b_0,\dots,b_{n-1},x_n)\Phi_1(dc|y_0,0,(x_1,b_1,\dots,x_n,\Delta,\Delta,\dots))~\mbox{on~} (0,\infty];\\
&&\psi_n(db|\hat{h}_n,c):=\frac{\pi_n^{u_0}(db|x_0,b_0,x_1,b_1,\dots,x_n)}{1-\pi^{u_0}_n(\{\Delta\}|x_0,b_0,x_1,b_1,\dots,x_n)}I\{c=0\}\\
&&+I\{c>0\}
\frac{\pi^{\Gamma^1_1,(y_0,0,(x_0,b_0,\dots,x_n,\Delta,\dots))}_0(db|x_n)}{1-\pi^{\Gamma^1_1,(y_0,0,(x_0,b_0,\dots,x_n,\Delta,\dots))}_0(\{\Delta\}|x_n)}\\
&&=\frac{\pi_n^{u_0}(db|x_0,b_0,x_1,b_1,\dots,x_n)}{1-\pi^{u_0}_n(\{\Delta\}|x_0,b_0,x_1,b_1,\dots,x_n)}I\{c=0\}+I\{c>0\}
\pi^{\Gamma^1_1,(y_0,0,(x_0,b_0,\dots,x_n,\Delta,\dots))}_0(db|x_n) ;\\
&&F^n(\hat{h}_n)_t(da):=\Pi_1(da|y_0,0,(x_0,b_0,\dots,x_n,\Delta,\Delta,\dots),t).
\end{eqnarray*}
Recall that $y_0=(x_0,\Delta,\Delta,\dots).$

If $\{1\le i\le n: ~\theta_i>0\}\ne \emptyset,$ then let $m(\hat{h}_n):=\#\{1\le i\le n: ~\theta_i>0\}$, and $l(\hat{h}_n):=\max\{1\le i\le n: \theta_i>0\}.$
When the context is clear, we write $m$ and $l$ instead of $m(\hat{h}_n)$ and $l(\hat{h}_n)$ for brevity.
Let $h_m$ be the $m$-history in the gradual-impulse control problem contained in $\hat{h}_n$. More precisely, $h_m$ is defined based on $\hat{h}_n$ as follows. Let $\tau_0(\hat{h}_n)=0,$ and $\tau_i(\hat{h}_n):=\inf\{j>\tau_{i-1}:~\theta_j>0\}$ for each $i\ge 1.$ Note that $l=\tau_m.$ Then $h_m=h_m(\hat{h}_n)=(y_0,0,y_1,\theta_{\tau_1},y_{2},\dots,\theta_{\tau_{m-1}},y_m),$ where
\begin{eqnarray*}
&&y_0=(x_0,\Delta,\Delta,\dots);~y_1=(x_0,b_0,x_1,b_1,\dots,x_{\tau_1-1},\Delta,\Delta,\dots);\\
&&\mbox{if $\theta_{\tau_1}=c_{\tau_1-1}$, then~} y_2=(x_{\tau_1-1},b_{\tau_1-1},x_{\tau_1},b_{\tau_1},\dots,x_{\tau_2-1},\Delta,\Delta,\dots),\\
&&\mbox{if $\theta_{\tau_1}<c_{\tau_1-1}$, then~} y_2=(x_{\tau_1},b_{\tau_1},\dots,x_{\tau_2-1},\Delta,\Delta,\dots);\\
&&\vdots\\
&&\mbox{if $\theta_{\tau_{m-1}}=c_{\tau_{m-1}-1}$, then~} y_m=(x_{\tau_{m-1}-1},\dots,x_{\tau_{m}-1},\Delta,\Delta,\dots),\\
&&\mbox{if $\theta_{\tau_{m-1}}<c_{\tau_{m-1}-1}$, then~} y_m=(x_{\tau_{m-1}},\dots,x_{\tau_{m}-1},\Delta,\Delta,\dots).
\end{eqnarray*}
For example, if
\begin{eqnarray*}\hat{h}_5=((0,x_0),(b_0,0,\rho^0),(0,x_1),(b_1,3,\rho^1),(3,x_2),(b_2,0,\rho^2),(0,x_3),(b_3,2,\rho^3),(1,x_4),(b_4,0,\rho^4),(0,x_5)),
 \end{eqnarray*}
 then $n=5$, $m=2$, $l=4,$ $\tau_1=2$, $\tau_2=4,$ and $h_2=(y_0,0,y_1,3,y_2)$ with
$y_1=(x_0,b_0,x_1,\Delta,\dots)$ and $y_2=(x_1,b_1,x_2,b_2,x_3,\Delta,\dots).$
Roughly speaking, the integer $m(\hat{h}_n)$ counts the number of interventions (except $y_0$) contained in the $n$-history of the hat DTMDP model.

 If $0<\theta_l=c_{l-1},$  we define
\begin{eqnarray*}
&&\varphi_n(\{0\}|\hat{h}_n):=1-\pi^{\Gamma^1_m,h_m}_{n-l+1}(\{\Delta\}|x_{l-1},b_{l-1},\dots,b_{n-1},x_n),\\
&&\varphi_n(dc|\hat{h}_n):=\pi^{\Gamma^1_m,h_m}_{n-l+1}(\{\Delta\}|x_{l-1},b_{l-1},\dots,b_{n-1},x_n)\Phi_m(dc|h_m)~\mbox{on~} (0,\infty];\\
&&\psi_n(db|\hat{h}_n,c):=\frac{\pi^{\Gamma^1_m,h_m}_{n-l+1}(db|x_{l-1},b_{l-1},\dots,b_{n-1},x_n)}{1-\pi^{\Gamma^1_m,h_m}_{n-l+1}(\{\Delta\}|x_{l-1},b_{l-1},\dots,b_{n-1},x_n)}I\{c=0\}\\
&&
+I\{c>0\}\frac{\pi^{\Gamma^1_{m+1},(h_m,\theta_l,(x_{l-1},b_{l-1},\dots,x_n,\Delta,\dots))}_0(db|x_n)}{1-\pi^{\Gamma^1_{m+1},(h_m,\theta_l,(x_{l-1},b_{l-1},\dots,x_n,\Delta,\dots))}_0(\{\Delta\}|x_n)}\\
&&= \frac{\pi^{\Gamma^1_m,h_m}_{n-l+1}(db|x_{l-1},b_{l-1},\dots,b_{n-1},x_n)}{1-\pi^{\Gamma^1_m,h_m}_{n-l+1}(\{\Delta\}|x_{l-1},b_{l-1},\dots,b_{n-1},x_n)}I\{c=0\}\\
&&+I\{c>0\} \pi^{\Gamma^1_{m+1},(h_m,\theta_l,(x_{l-1},b_{l-1},\dots,x_n,\Delta,\dots))}_0(db|x_n);\\
&&F^n(\hat{h}_n)_t(da):=\Pi_m(da|h_m,t).
\end{eqnarray*}
Finally, if $0<\theta_l<c_{l-1}$, then we define
\begin{eqnarray*}
&&\varphi_n(\{0\}|\hat{h}_n):=1-\pi^{\Gamma^0_m,h_m,\theta_l,x_l}_{n-l}(\{\Delta\}|x_{l},b_{l},\dots,b_{n-1},x_n),\\
&&\varphi_n(dc|\hat{h}_n):=\pi^{\Gamma^0_m,h_m,\theta_l,x_l}_{n-l}(\{\Delta\}|x_{l},b_{l},\dots,b_{n-1},x_n)\Phi_m(dc|h_m)~\mbox{on~} (0,\infty];\\
&&\psi_n(db|\hat{h}_n,c):=\frac{\pi^{\Gamma^0_m,h_m,\theta_l,x_l}_{n-l}(db|x_{l},b_{l},\dots,b_{n-1},x_n)}{1-\pi^{\Gamma^0_m,h_m,\theta_l,x_l}_{n-l}(\{\Delta\}|x_{l},b_{l},\dots,b_{n-1},x_n)}I\{c=0\}\\
&&
+I\{c>0\}\frac{\pi^{\Gamma^1_{m+1},(h_m,\theta_l,(x_{l},b_{l},\dots,x_n,\Delta,\dots))}_0(db|x_n)}{1-\pi^{\Gamma^1_{m+1},(h_m,\theta_l,(x_{l},b_{l},\dots,x_n,\Delta,\dots))}_0(\{\Delta\}|x_n)}\\
&=&\frac{\pi^{\Gamma^0_m,h_m,\theta_l,x_l}_{n-l}(db|x_{l},b_{l},\dots,b_{n-1},x_n)}{1-\pi^{\Gamma^0_m,h_m,\theta_l,x_l}_{n-l}(\{\Delta\}|x_{l},b_{l},\dots,b_{n-1},x_n)}I\{c=0\}\\
&&
+I\{c>0\} \pi^{\Gamma^1_{m+1},(h_m,\theta_l,(x_{l},b_{l},\dots,x_n,\Delta,\dots))}_0(db|x_n);\\
&&F^n(\hat{h}_n)_t(da):=\Pi_m(da|h_m,t).
\end{eqnarray*}

To be specific, we call the (typical) strategy $\sigma=\{(\varphi_n,\psi_n,F^n)\}_{n=0}^\infty$ defined above as the strategy induced by the policy $u.$ The next statement reveals a connection between a policy $u$ and its induced strategy $\sigma$ for the hat DTMDP model.
\begin{proposition}\label{JapanProposition01}
For each policy $u$ and the strategy $\sigma=\{(\varphi_n,\psi_n,F^n)\}_{n=0}^\infty$ induced by $u$,
${\cal V}(x,u)=V((0,x),\sigma)$, and therefore, ${\cal V}^\ast(x)\ge V^\ast(x)$ for each $x\in\textbf{X}.$
\end{proposition}

%\par\noindent\textit{Proof.} One can verify \begin{eqnarray*}
%&&{\rm E}_x^u\left[e^{\sum_{i=1}^n C^I(Y_i) +\sum_{i=2}^n\int_{0}^{\Theta_{i}} \int_{\textbf{A}^G}c^G(\overline{x}(Y_{i-1}),a)\Pi_{i-1}(da|H_{i-1},s)ds}\right]\\
%&=&\mathbb{E}_{(0,x)}^{\sigma}\left[e^{\sum_{i=0}^{\tau_n-1} c^I(X_i,B_i,X_{i+1}) +\sum_{i=2}^{n}\int_{0}^{\Theta_{\tau_{i-1}}} \int_{\textbf{A}^G}c^G(X_{\tau_{i-1}-1},a)F^{\tau_{i-1}-1}(\hat{H}_{\tau_{i-1}-1})_s(da)ds}\right]
%\end{eqnarray*}
%for each $n\ge 1$. The case of $n=1$ can be readily seen (we accept $\sum_{n=2}^1(\cdot):=0$), as a consequence of the definitions of the strategy $\sigma=\{(\varphi_n,\psi_n,F^n)\}_{n=0}^\infty$ induced by $u$. The general case follows from an inductive argument. The cumbersome details are omitted. Passing to the limit as $n\rightarrow \infty$ and an application of the monotone convergence theorem yield the equality in the statement. The last assertion holds automatically from the first assertion. $\hfill\Box$
%\bigskip

\begin{remark}\label{JapanRevisionAARemark01}
 A deterministic stationary policy say $u^D$ is associated with a strategy $\sigma^D=(\varphi,\psi,F)$ in the hat DTMDP model, where $F(x)_t(da)=\delta_{f(x)}(da)$ for all $t\ge 0.$ It is evident that ${\cal V}(x,u^D)=V(x,\sigma^D)$ for each $x\in\textbf{X}.$
Thus, if the hat DTMDP problem (\ref{JapanProblem02}) has an optimal strategy in the form of $\sigma^D=(\varphi,\psi,F)$, then the previous discussions lead to ${\cal V}^\ast(x)=V^\ast(x)$, and that the deterministic stationary policy $u^D$ associated with $\sigma^D$ is optimal for the gradual-impulse control problem (\ref{JapanProblem01}).
\end{remark}

\section{Proof of the main statements}\label{JapanSecProof}
In this section, we prove the results stated in Section \ref{JapanSecMain}. This is based on the investigation of problem (\ref{JapanProblem02}) for the hat DTMDP model described in Section \ref{JapanSecDTMDPhat}. In this section, unless specified otherwise, $V^\ast$ is understood as the value function of problem (\ref{JapanProblem02}) for the hat DTMDP model. More exactly, the main properties concerning $V^\ast$ are summarized in the next statement.
%\begin{remark}
%{We suppose that Condition \ref{JapanCon02} is always satisfied in this section, without special reference.}
%\end{remark}

\begin{proposition}\label{JapanTheorem01}
\begin{itemize}
\item[(a)] $V^\ast$ is a $[1,\infty]$-valued lower semianalytic function on $\textbf{X}$ satisfying
\begin{eqnarray}\label{JapanEqnA}
&&\inf_{a\in\textbf{A}^G}\left\{\int_{\textbf{X}} V^\ast(y) \tilde{q}(dy|x,a)-(q_x(a)-c^G(x,a))V^\ast(x)  \right\}\ge 0,\\
&&~\forall~x\in\textbf{X}^\ast(V^\ast):=\{x\in\textbf{X}:~V^\ast(x)<\infty\}\nonumber
\end{eqnarray}
and
\begin{eqnarray}\label{JapanEqnB}
V^\ast(x)\le \inf_{b\in\textbf{A}^I}\left\{ \int_{\textbf{X}} e^{c^I(x,b,y)}V^\ast(y)Q(dy|x,b) \right\},~x\in\textbf{X},
\end{eqnarray}
whereas at each $x\in\textbf{X}$,  the inequality in either (\ref{JapanEqnA}) or (\ref{JapanEqnB}) holds with equality.
\item[(b)] $\textbf{X}\setminus \textbf{X}^I\subseteq \textbf{X}^G,$ where $\textbf{X}^G:=\textbf{X}^G(V^\ast)$, see (\ref{JapanRevisionAAEqn02}), and $\textbf{X}^I:=\textbf{X}^I(V^\ast)$. ({Lemma \ref{JapanNewLem01} below asserts that $V^\ast$ is universally measurable so that the integrals $\int_{\textbf{X}}V^\ast(y)\tilde{q}(dy|x,a)$ and $\int_{\textbf{X}}V^\ast(y)e^{c^I(x,b,y)}Q(dy|x,b)$ are defined.})
%\item[(c)] There exist measurable mappings $\psi^\ast$ and $f^\ast$ from $\textbf{X}$ to $\textbf{A}^I$ and $\textbf{A}^G$, respectively, such that
%\begin{eqnarray*}
%&&\inf_{a\in \textbf{A}^G}\left\{\int_{\textbf{X}}V^\ast(y)\tilde{q}(dy|x,a)-(q_x(a)-c^G(x,a))V^\ast(x)\right\}\\
%&=&\int_{\textbf{X}}V^\ast(y)\tilde{q}(dy|x,f^\ast(x))-(q_x(f^\ast(x))-c^G(x,f^\ast(x)))V^\ast(x)
%\end{eqnarray*}
%for each $x\in\textbf{X}^G$, and
%\begin{eqnarray*}
%\inf_{b\in\textbf{A}^I}\left\{\int_{\textbf{X}}e^{c^I(x,b,y)}V^\ast(y)Q(dy|x,b)\right\}=\int_{\textbf{X}}V^\ast(y)e^{c^I(x,\psi^\ast(x),y)}Q(dy|x,\psi^\ast(x)),~\forall~x\in\textbf{X}.
%\end{eqnarray*}
%\item[(d)] For each pair of measurable mappings $(\psi^\ast,f^\ast)$ that satisfy the previous two relations,
\end{itemize}
\end{proposition}
\par\noindent\textit{Proof.} See Lemmas \ref{JapanNewLem01}, \ref{JapanLemma07} and \ref{JapanLemma06} below. $\hfill\Box$
\bigskip

%We introduce the notation to be used in the next statement:
%\begin{eqnarray*}
%\textbf{X}^G:=\left\{x\in\textbf{X}:\infty>V^\ast(x)=\inf_{a\in \textbf{A}^G}\left\{\int_{\textbf{X}}V^\ast(y)\tilde{q}(dy|x,a)-(q_x(a)-c^G(x,a))V^\ast(x)\right\}\right\},
%\end{eqnarray*}
%and denote by $\textbf{X}^I$ the collection of $x\in\textbf{X}$ at which, %$V^\ast(x)=\inf_{b\in\textbf{A}^I}\left\{\int_{\textbf{X}}V^\ast(y)e^{c^I(x,b,y)}Q(dy|x,b)\right\}.$

\begin{lemma}\label{JapanNewLem01}
The following assertions hold.
\begin{itemize}
\item[(a)] The value function $V^\ast$ depends on the state $(\theta,x)$ only through the second coordinate, and thus we write $V^\ast(x)$ instead of $V^\ast(\theta,x).$ The function $V^\ast$ is an $[1,\infty]$-valued lower semianalytic function satisfying
\begin{eqnarray}\label{JapanBellman01}
V(x)&=&\inf_{\hat{a}\in\hat{\textbf{A}}}\left\{\int_0^c \int_{\textbf{X}}V(y)\tilde{q}(dy|x,\rho_t)e^{-\int_0^t (q_x(\rho_s)-c^G(x,\rho_s))ds}dt\right. \\
&&\left. +I\{c=\infty\} e^{-\int_0^\infty q_x(\rho_s)ds} e^{\int_0^\infty c^G(x,\rho_s)ds}\right.\nonumber\\
&&\left.+I\{c<\infty\}e^{-\int_0^c (q_x(\rho_s)-c^G(x,\rho_s))ds}\int_{\textbf{X}}V(y)e^{c^I(x,b,y)}Q(dy|x,b)\right\},~x\in\textbf{X};\nonumber\\
V(x_\infty)&=&1,\nonumber
\end{eqnarray}
and is the minimal $[1,\infty]$-valued lower semianalytic function satisfying the following inequality
 \begin{eqnarray}\label{JapanBellman02}
V(x)&\ge&\inf_{\hat{a}\in\hat{\textbf{A}}}\left\{\int_0^c \int_{\textbf{X}}V(y)\tilde{q}(dy|x,\rho_t)e^{-\int_0^t (q_x(\rho_s)-c^G(x,\rho_s))ds}dt \right.\\
&&\left.+I\{c=\infty\} e^{-\int_0^\infty q_x(\rho_s)ds} e^{\int_0^\infty c^G(x,\rho_s)ds}\right.\nonumber\\
&&\left.+I\{c<\infty\}e^{-\int_0^c (q_x(\rho_s)-c^G(x,\rho_s))ds}\int_{\textbf{X}}V(y)e^{c^I(x,b,y)}Q(dy|x,b)\right\},~x\in\textbf{X};\nonumber\\
V(x_\infty)&=&1.\nonumber
\end{eqnarray}
\item[(b)] For each $\epsilon>0,$ there exists an $\epsilon$-optimal deterministic Markov universally measurable strategy that depends on the state $(\theta,x)$ only through the second coordinate for the hat DTMDP problem (\ref{JapanProblem02}). (The meaning of universally measurable strategies can be found in Appendix \ref{JapanAppendixA}.)
\item[(c)] A deterministic stationary strategy that depends on the state $(\theta,x)$ only through $x$ is optimal if and only if it attains the infimum in (\ref{JapanBellman01}) with $V^\ast$ replacing $V$, for each $x\in\textbf{X}.$
\item [(d)] For each $x\in\textbf{X}$, $V^\ast(x)=\inf_{\pi\in\Pi^{U}}V(x,\pi)$, where $\Pi^U$ indicates the class of universally measurable strategies in the hat DTMDP model.
\end{itemize}
\end{lemma}
\par\noindent\textit{Proof.} The fact that the value function $V^\ast$ is the minimal $[1,\infty]$-valued lower semianalytic function satisfying \begin{eqnarray*}
g(\theta,x)&\ge&\inf_{\hat{a}\in\hat{\textbf{A}}}\left\{\int_0^c \int_{\textbf{X}}g(t,y)\tilde{q}(dy|x,\rho_t)e^{-\int_0^t (q_x(\rho_s)-c^G(x,\rho_s))ds}dt\right.\\
&&\left. +I\{c=\infty\} e^{-\int_0^\infty q_x(\rho_s)ds} e^{\int_0^\infty c^G(x,\rho_s)ds}\right.\\
&&\left.+I\{c<\infty\}e^{-\int_0^c (q_x(\rho_s)-c^G(x,\rho_s))ds}\int_{\textbf{X}}g(c,y)e^{c^I(x,by)}Q(dy|x,b)\right\},~x\in\textbf{X};\\
g(\infty,x_\infty)&:=&1,
\end{eqnarray*}
where the inequality can be replaced by equality,
follows from Proposition \ref{GGZyExponentialProposition02}. The existence of an $\epsilon$-optimal deterministic Markov universally measurable strategy follows from Proposition \ref{GGZyExponentialProposition02}, too. Furthermore,
note that the first coordinate in the state $(\theta,x)$ does not affect the cost function or the transition probability, from which the independence on the first coordinate of the state $(\theta,x)$ follows, c.f. \cite{Feinberg:2005}. Now assertions (a,b) follow. Finally, the last two assertions follow from Proposition \ref{GGZyExponentialProposition02}. $\hfill\Box$
\bigskip

\begin{lemma}\label{JapanLemma02}
The function
\begin{eqnarray*}
t\in[0,\infty)\rightarrow \int_0^t \int_{\textbf{X}} e^{-\int_0^\tau (q_x(\rho_s)-c^G(x,\rho_s))ds}V^\ast(y)\tilde{q}(dy|x,\rho_\tau)d\tau +e^{-\int_0^t (q_x(\rho_s)-c^G(x,\rho_s))ds}V^\ast(x)
\end{eqnarray*}
\end{lemma}
is increasing, for each $x\in\textbf{X}$ and $\rho\in{\cal R}.$

\par\noindent\textit{Proof.}
Let $0\le t_1<t_2<\infty$ and $x\in\textbf{X}$ be fixed, and we will verify
\begin{eqnarray*}
&&\int_0^{t_2} e^{-\int_0^\tau (q_x(\rho_s)-c^G(x,\rho_s))ds}\int_{\textbf{X}}V^\ast(y)\tilde{q}(dy|x,\rho_\tau)d\tau+e^{-\int_0^{t_2} (q_x(\rho_s)-c^G(x,\rho_s))ds}V^\ast(x)\\
&\ge&\int_0^{t_1} e^{-\int_0^\tau (q_x(\rho_s)-c^G(x,\rho_s))ds}\int_{\textbf{X}}V^\ast(y)\tilde{q}(dy|x,\rho_\tau)d\tau+e^{-\int_0^{t_1} (q_x(\rho_s)-c^G(x,\rho_s))ds}V^\ast(x),
\end{eqnarray*}
as follows. It is sufficient to consider the case when the left hand side is finite, for otherwise, the above inequality would hold automatically. Then the goal is to show, by subtracting the right hand side from the left hand side,
\begin{eqnarray*}
&&0\le \int_{t_1}^{t_2} e^{-\int_0^\tau (q_x(\rho_s)-c^G(x,\rho_s))ds}\int_{\textbf{X}}V^\ast(y)\tilde{q}(dy|x,\rho_\tau)d\tau + e^{-\int_0^{t_2} (q_x(\rho_s)-c^G(x,\rho_s))ds}V^\ast(x)\\
&&-e^{-\int_0^{t_1} (q_x(\rho_s)-c^G(x,\rho_s))ds}V^\ast(x).
\end{eqnarray*}
The right hand side of this inequality can be further written as
\begin{eqnarray*}
&&\int_0^{t_2-t_1} e^{-\int_0^{t_1} (q_x(\rho_s)-c^G(x,\rho_s))ds  }e^{-\int_{t_1}^{\tau+t_1}(q_x(\rho_s)-c^G(x,\rho_s))ds }\int_{\textbf{X}} V^\ast(y)\tilde{q}(dy|x,\rho_{\tau+t_1})d\tau \\
&&+e^{-\int_{0}^{t_1}(q_x(\rho_s)-c^G(x,\rho_s))ds }\left(e^{-\int_{t_1}^{t_2}(q_x(\rho_s)-c^G(x,\rho_s))ds}-1\right)V^\ast(x)\\
&=&e^{-\int_{0}^{t_1}(q_x(\rho_s)-c^G(x,\rho_s))ds } \left\{\int_0^{t_2-t_1} e^{-\int_0^\tau (q_x(\rho_{s+t_1})-c^G(x,\rho_{s+t_1}))ds}\int_{\textbf{X}}V^\ast(y)\tilde{q}(dy|x,\rho_{\tau+t_1})d\tau \right.\\
&&\left.+\left(e^{-\int_{0}^{t_2-t_1 }(q_x(\rho_{t_1+s})-c^G(x,\rho_{t_1+s}))ds}-1\right)V^\ast(x) \right\}.
\end{eqnarray*}
Introduce $\tilde{\rho}_s:=\rho_{t_1+s}$ for each $s\ge 0.$ The target becomes to show
\begin{eqnarray*}
\int_0^{t_2-t_1} e^{-\int_0^\tau (q_x(\tilde{\rho}_s)-c^G(x,\tilde{\rho}_s))ds}\int_{\textbf{X}}V^\ast(y)\tilde{q}(dy|x,\tilde{\rho}_\tau)d\tau +e^{-\int_0^{t_2-t_1}(q_x(\tilde{\rho}_s)-c^G(x,\tilde{\rho}_s))ds }V^\ast(x)\ge V^\ast(x).
\end{eqnarray*}
To this end, for a fixed $\epsilon>0$, let us consider a deterministic Markov $\epsilon$-optimal universally measurable strategy $\{(\varphi^\ast_n,\psi^\ast_n,F^{\ast,n})\}_{n=0}^\infty$ coming from Lemma \ref{JapanNewLem01}, and an associated universally measurable strategy $\pi^{New}=\{(\varphi_n,\psi_n,F^n)\}_{n=0}^\infty$ defined by
$\varphi_0(\theta, x):=\varphi_0^\ast(x)+t_2-t_1,$ $\psi_0(\theta,x)=\psi^\ast_0(x)$, $F^0(\theta,x)_s=\tilde{\rho}_s$ if $s\le t_2-t_1$ and $F^0(\theta,x)_s=F^{\ast,0}(\theta,x)_{s-(t_2-t_1)}$ if $s>t_2-t_1;$ and for $n\ge 1$, $\varphi_n((\theta,x),\hat{a},(t,y))= \varphi^\ast_{n-1}(y)$, $\psi_n((\theta,x),\hat{a},(t,y))=\psi^\ast_{n-1}(y),$ and $F^n((\theta,x),\hat{a},(t,y))_s= F^{\ast,n-1}(y)_s$ for all $s\ge 0.$ Under the universally measurable strategy $\pi^{New}$, only the gradual control action $\tilde{\rho}$ is used up to either $t_2-t_1$ or the natural jump moment, whichever takes place first, after when, the $\epsilon$-optimal universally measurable strategy is in use, and so
\begin{eqnarray*}
&&V^\ast(x)\le V(x,\pi^{New})\\
&\le& \int_0^{t_2-t_1} e^{-\int_0^\tau (q_x(\tilde{\rho}_s) -c^G(x,\tilde{\rho}_s))ds}\int_{\textbf{X}}(V^\ast(y)+\epsilon)\tilde{q}(dy|x,\tilde{\rho}_\tau)d\tau+e^{-\int_0^{t_2-t_1}(q_x(\tilde{\rho}_s)-c^{G}(x,\tilde{\rho}_s))ds }(V^\ast(x)+\epsilon)\\
&=&\int_0^{t_2-t_1} e^{-\int_0^\tau (q_x(\tilde{\rho}_s) -c^G(x,\tilde{\rho}_s))ds}\int_{\textbf{X}}V^\ast(y)\tilde{q}(dy|x,\tilde{\rho}_\tau)d\tau+e^{-\int_0^{t_2-t_1}(q_x(\tilde{\rho}_s)-c^{G}(x,\tilde{\rho}_s))ds }V^\ast(x)\\
&&+\epsilon\left( \int_0^{t_2-t_1} e^{-\int_0^\tau (q_x(\tilde{\rho}_s) -c^G(x,\tilde{\rho}_s))ds}q_x(\tilde{\rho}_\tau) )d\tau+e^{-\int_0^{t_2-t_1}(q_x(\tilde{\rho}_s)-c^{G}(x,\tilde{\rho}_s))ds}\right),
\end{eqnarray*}
where the first inequality holds because of the last assertion of Lemma \ref{JapanNewLem01}.
Since the expression in the last bracket is nonnegative and finite, and $\epsilon>0$ was arbitrarily fixed,
we see that $V^\ast(x) \le \int_0^{t_2-t_1} e^{-\int_0^\tau (q_x(\tilde{\rho}_s) -c^G(x,\tilde{\rho}_s))ds}\int_{\textbf{X}}V^\ast(y)\tilde{q}(dy|x,\tilde{\rho}_\tau)d\tau+e^{-\int_0^{t_2-t_1}(q_x(\tilde{\rho}_s)-c^{G}(x,\tilde{\rho}_s))ds }V^\ast(x),$ as desired. $\hfill\Box$
\bigskip

\begin{lemma}\label{JapanLemma07}
Relations (\ref{JapanEqnA}) and (\ref{JapanEqnB}) hold. (Recall from Lemma \ref{JapanNewLem01} that $V^\ast$ is universally measurable.)
\end{lemma}
\par\noindent\textit{Proof.}  Let $x\in\textbf{X}$ be fixed.
Inequality (\ref{JapanEqnB}) immediately follows from Lemma \ref{JapanNewLem01}, if on the right hand side of (\ref{JapanBellman01}) with $V^\ast$ replacing $V$, one takes the infimum over actions $\hat{a}\in\hat{\textbf{A}}$ with $c=0$. (Recall the notation in use: $\hat{a}=(c,b,\rho)\in\hat{\textbf{A}}$.) Let us verify (\ref{JapanEqnA}) as follows.  Suppose $V^\ast(x)<\infty$. Let $a\in\textbf{A}^G$ be arbitrarily fixed. If $\int_{\textbf{X}}V^\ast(y)\tilde{q}(dy|x,a)=\infty$, then trivially,
$\int_{\textbf{X}} V^\ast(y) \tilde{q}(dy|x,a)-(q_x(a)-c^G(x,a))V^\ast(x)\ge 0$. Consider the case when $\int_{\textbf{X}}V^\ast(y)\tilde{q}(dy|x,a)<\infty.$ Let $t>0$ be arbitrarily fixed. Then $\int_0^t e^{-\tau (q_x(a) -c^G(x,a))} \int_{\textbf{X}}V^\ast(y)\tilde{q}(dy|x,a)d\tau+e^{-t (q_x(a)-c^G(x,a))}V^\ast(x)$ is finite. Upon differentiating it with respect to $t$ and applying the fundamental theorem of calculus, we see
 \begin{eqnarray*}
 e^{-(q_x(a)-c^G(x,a))t} \int_{\textbf{X}} V^\ast(y)\tilde{q}(dy|x,a)-(q_x(a)-c^G(x,a))e^{-t(q_x(a)-c^G(x,a))}V^\ast(x)\ge 0,
 \end{eqnarray*}
where the inequality follows from Lemma \ref{JapanLemma02}. Thus, $\int_{\textbf{X}} V^\ast(y) \tilde{q}(dy|x,a)-(q_x(a)-c^G(x,a))V^\ast(x)\ge 0$. Since $a\in\textbf{A}^G$ was arbitrarily fixed, we see that (\ref{JapanEqnA}) holds.
$\hfill\Box$
\bigskip

%(Under the assumptions of the current statement, if $V^\ast(x)=\infty,$ then necessarily for all $a\in\textbf{A}^G$,  $\int_{\textbf{X}}V^\ast(y)\tilde{q}(dy|x,a)=\infty$. Indeed, suppose for contradiction that for some $a\in\textbf{A}^G$, $\int_{\textbf{X}}V^\ast(y)\tilde{q}(dy|x,a)<\infty$. Then the strategy that identically applies the gradual control $a$ and does not apply any impulsive control would be with a finite value, which is against $V^\ast(x)=\infty$. Thus, with the convenition of $\infty-\infty:=\infty$, (\ref{JapanEqnA}) and (\ref{JapanEqnB}) hold for all $x\in\textbf{X}.$)

\begin{lemma}\label{JapanLemma06}
For each $x\in \textbf{X}$, the inequality in either (\ref{JapanEqnA}) or (\ref{JapanEqnB}) holds with equality.
\end{lemma}
\par\noindent\textit{Proof.} Let $x\in\textbf{X}$ be fixed. If the equality in (\ref{JapanEqnB}) holds at this point, then there is nothing to prove. Suppose the strict inequality holds in (\ref{JapanEqnB}). Then necessarily $V^\ast(x)<\infty,$ and so $x\in {\bf X}^\ast(V^\ast).$ The objective is to show that, in this case, (\ref{JapanEqnA}) holds with equality. For the infimum in (\ref{JapanBellman01}) with $V^\ast$ replacing $V$, it suffices to consider $c>0$, because (\ref{JapanEqnB}) holds with strict inequality at the fixed point $x\in\textbf{X}$ here.

{Suppose for contradiction that the strict inequality holds in (\ref{JapanEqnA}), and we let
\begin{eqnarray*}
\Xi:=\inf_{a\in\textbf{A}^G}\left\{\int_{\textbf{X}} V^\ast(y) \tilde{q}(dy|x,a)-(q_x(a)-c^G(x,a))V^\ast(x)  \right\}> 0.
\end{eqnarray*}}

{Since by assumption  the strict inequality holds in (\ref{JapanEqnB}) at the given point $x\in {\bf X}$,
 there is some $\delta>0$ such that
\begin{eqnarray*}
&&\delta \bar{q}_x<\Xi;\\
&& \inf_{b\in\textbf{A}^I}\left\{ \int_{\textbf{X}} e^{c^I(x,b,y)}V^\ast(y)Q(dy|x,b) \right\}-V^\ast(x)>\delta.
\end{eqnarray*}
}

{Let $T\in(0,\infty)$ be fixed.} Let {
\begin{eqnarray*}
\epsilon\in \left(0,\min\left\{\delta,\int_0^T e^{-\bar{q}_x s} ds \Xi\right\}\right)
\end{eqnarray*}}
be fixed, and $(c^\ast,b^\ast,\rho^\ast)\in\hat{\textbf{A}}$ be such that
\begin{eqnarray*}
&&V^\ast(x)+\epsilon\\
&\ge &\left\{\int_0^{c^\ast} \int_{\textbf{X}}V^\ast(y)\tilde{q}(dy|x,\rho^\ast_t)e^{-\int_0^t (q_x(\rho^\ast_s)-c^G(x,\rho^\ast_s))ds}dt +I\{c^\ast=\infty\} e^{-\int_0^\infty q_x(\rho^\ast_s)ds} e^{\int_0^\infty c^G(x,\rho^\ast_s)ds}\right.\nonumber\\
&&\left.+I\{c^\ast<\infty\}e^{-\int_0^{c^\ast} (q_x(\rho^\ast_s)-c^G(x,\rho^\ast_s))ds}\int_{\textbf{X}}V^\ast(y)e^{c^I(x,b^\ast,y)}Q(dy|x,b^\ast)\right\}.
\end{eqnarray*}
There are two cases to be considered: (a) $0<c^\ast<\infty$; (b) $c^\ast=\infty$.

Consider case (a).

 {Since $V^\ast(x)<\infty$, \begin{eqnarray*}\int_0^{c^\ast} e^{-{\int_0^t (q_x(\rho^\ast_s)-c^G(x,\rho^\ast_s))ds }} \int_{\textbf{X}} V^\ast(y)\tilde{q}(dy|x,\rho^\ast_t)dt<\infty,~-\infty<\int_0^{c^\ast} (q_x(\rho^\ast_s)-c^G(x,\rho^\ast_s))ds<\infty,
\end{eqnarray*}
and
\begin{eqnarray*}
\epsilon+V^\ast(x)&\ge &\int_0^{c^\ast} \int_{\textbf{X}}V^\ast(y)\tilde{q}(dy|x,\rho^\ast_t)e^{-\int_0^t (q_x(\rho^\ast_s)-c^G(x,\rho^\ast_s))ds}dt\\
&&+
e^{-\int_0^{c^\ast} (q_x(\rho^\ast_s)-c^G(x,\rho^\ast_s))ds}\int_{\textbf{X}}V^\ast(y)e^{c^I(x,b^\ast,y)}Q(dy|x,b^\ast)\\
&\ge&\int_0^{c^\ast} \int_{\textbf{X}}V^\ast(y)\tilde{q}(dy|x,\rho^\ast_t)e^{-\int_0^t (q_x(\rho^\ast_s)-c^G(x,\rho^\ast_s))ds}dt\\
&&+
e^{-\int_0^{c^\ast} (q_x(\rho^\ast_s)-c^G(x,\rho^\ast_s))ds}(V^\ast(x)+\delta) \\
&=&\int_0^{c^\ast} \int_{\textbf{X}}V^\ast(y)\tilde{q}(dy|x,\rho^\ast_t)e^{-\int_0^t (q_x(\rho^\ast_s)-c^G(x,\rho^\ast_s))ds}dt\\
&&+
e^{-\int_0^{c^\ast} (q_x(\rho^\ast_s)-c^G(x,\rho^\ast_s))ds}V^\ast(x)+\delta e^{-\int_0^{c^\ast} (q_x(\rho^\ast_s)-c^G(x,\rho^\ast_s))ds} \\
&=&\int_0^{c^\ast} \left\{\int_{\textbf{X}}V^\ast(y)\tilde{q}(dy|x,\rho^\ast_t)-(q_x(\rho^\ast_t)-c^G(x,\rho^\ast_t))V^\ast(x) \right\} e^{-\int_0^t (q_x(\rho^\ast_s)-c^G(x,\rho^\ast_s))ds}dt\\
&&+V^\ast(x)+\delta e^{-\int_0^{c^\ast} (q_x(\rho^\ast_s)-c^G(x,\rho^\ast_s))ds}\\
&\ge& \int_0^{c^\ast} \Xi e^{-\int_0^t (q_x(\rho^\ast_s)-c^G(x,\rho^\ast_s))ds}dt\\
&&+V^\ast(x)+\delta -\int_0^{c^\ast} e^{-\int_0^{t} (q_x(\rho^\ast_s)-c^G(x,\rho^\ast_s))ds}(q_x(\rho^\ast_t)-c^G(x,\rho^\ast_t))\delta dt,
\end{eqnarray*}
i.e.,
\begin{eqnarray*}
\epsilon&\ge&\int_0^{c^\ast} \left(\Xi-(q_x(\rho^\ast_t)-c^G(x,\rho^\ast_t))\delta\right) e^{-\int_0^t (q_x(\rho^\ast_s)-c^G(x,\rho^\ast_s))ds}dt+\delta\\
&\ge&\int_0^{c^\ast} \left(\Xi-q_x(\rho^\ast_t)\delta\right) e^{-\int_0^t (q_x(\rho^\ast_s)-c^G(x,\rho^\ast_s))ds}dt+\delta\\
&\ge& \delta,
\end{eqnarray*}
where the last inequality holds because $\Xi-q_x(\rho^\ast_t)\delta\ge \Xi-\bar{q}_x\delta>0$. This contradicts $\epsilon<\delta$, as desired.}

{Therefore, \begin{eqnarray*}
\inf_{a\in\textbf{A}^G}\left\{\int_{\textbf{X}} V^\ast(y) \tilde{q}(dy|x,a)-(q_x(a)-c^G(x,a))V^\ast(x)  \right\}=0
\end{eqnarray*}
in case (a).}

{Now consider case (b). Then
\begin{eqnarray*}
\epsilon+V^\ast(x)\ge \int_0^{\infty} e^{-{\int_0^t (q_x(\rho^\ast_s)-c^G(x,\rho^\ast_s))ds }} \int_{\textbf{X}} V^\ast(y)\tilde{q}(dy|x,\rho^\ast_t)dt+e^{-\int_0^{\infty} q_x(\rho^\ast_s)ds} e^{\int_0^\infty c^G(x,\rho^\ast_s)ds}.
\end{eqnarray*}}
{One can apply the proof of Lemma 5.3 of \cite{Zhang:2017} to showing that for each $t\in[0,\infty),$
\begin{eqnarray}\label{JapanEasyReadingEqn01}
V^\ast(x)+\epsilon\ge \inf_{\rho\in{\cal R}}\left\{\int_0^{t} e^{-{\int_0^t (q_x(\rho_s)-c^G(x,\rho_s))ds }} \int_{\textbf{X}} V^\ast(y)\tilde{q}(dy|x,\rho_t)dt+e^{-\int_0^{t} (q_x(\rho_s)-c^G(x,\rho_s))ds}V^\ast(x)   \right\}.
\end{eqnarray}
To improve the readability, we provide the detailed justification of this fact as follows. We only need consider when $t>0$; the case of $t=0$ is trivial. Let $\delta'>0$ be arbitrarily fixed. Then there is some $\hat{\rho}\in {\cal R}$ such that
\begin{eqnarray*}
\epsilon+V^\ast(x)+\delta'\ge \int_0^\infty e^{-\int_0^\tau (q_x(\hat{\rho}_s)-c^G(x,\hat{\rho}_s))ds}\int_{\textbf{X}}V^\ast(y)\tilde{q}(dy|x,\hat{\rho}_\tau)d\tau+e^{-\int_0^\infty q_x(\hat{\rho}_s)ds}e^{\int_0^\infty c^G(x,\hat{\rho}_s)ds}.
\end{eqnarray*}
Define $\tilde{\rho}\in {\cal R}$ by $\tilde{\rho}_s=\hat{\rho}_{t+s}$ for each $s\ge 0.$
Then, for each $t\ge 0,$
\begin{eqnarray*}
&& \epsilon+V^\ast(x)+\delta'\\
&\ge& \int_0^t e^{-\int_0^\tau (q_x(\hat{\rho}_s)-c^G(x,\hat{\rho}_s))ds}\int_{\textbf{X}}V^\ast(y)\tilde{q}(dy|x,\hat{\rho}_\tau)d\tau +\int_t^\infty e^{-\int_0^\tau (q_x(\hat{\rho}_s)-c^G(x,\hat{\rho}_s))ds}\int_{\textbf{X}}V^\ast(y)\tilde{q}(dy|x,\hat{\rho}_\tau)d\tau\\
&&+e^{-\int_0^t (q_x(\hat{\rho}_s)-c^G(x,\hat{\rho}_s))ds} e^{-\int_t^\infty q_x(\hat{\rho}_s)ds}e^{\int_t^\infty c^G(x,\hat{\rho}_s)ds}\\
&=& \int_0^t e^{-\int_0^\tau (q_x(\hat{\rho}_s)-c^G(x,\hat{\rho}_s))ds}\int_{\textbf{X}}V^\ast(y)\tilde{q}(dy|x,\hat{\rho}_\tau)d\tau +e^{-\int_0^t (q_x(\hat{\rho}_v)-c^G(x,\hat{\rho}_v))dv}\\
&&\times \left\{\int_0^\infty e^{-\int_0^s (q_x(\tilde{\rho}_v))-c^G(x,\tilde{\rho}_v))dv}\int_{\textbf{X}}V^\ast(y)\tilde{q}(dy|x,\tilde{\rho}_s)ds+ e^{-\int_0^\infty q_x(\tilde{\rho}_s)ds}e^{\int_0^\infty c^G(x,\tilde{\rho}_s)ds} \right\}\\
&\ge &\int_0^t e^{-\int_0^\tau (q_x(\hat{\rho}_s)-c^G(x,\hat{\rho}_s))ds}\int_{\textbf{X}}V^\ast(y)\tilde{q}(dy|x,\hat{\rho}_\tau)d\tau +e^{-\int_0^t (q_x(\hat{\rho}_v)-c^G(x,\hat{\rho}_v))dv}V^\ast(x)\\
&\ge& \inf_{\rho\in{\cal R}} \left\{\int_0^t e^{-\int_0^\tau (q_x(\rho_s)-c^G(x,{\rho}_s))ds}\int_{\textbf{X}}V^\ast(y)\tilde{q}(dy|x,{\rho}_\tau)d\tau +e^{-\int_0^t (q_x({\rho}_v)-c^G(x,{\rho}_v))dv}V^\ast(x)\right\},
\end{eqnarray*}
where the second inequality is by Lemma \ref{JapanNewLem01}(a), which in particular, asserts that $V^\ast$ satisfies (\ref{JapanBellman01}). Since $\delta'>0$ was arbitrarily fixed, the above implies (\ref{JapanEasyReadingEqn01}).}

{Fix $\delta''>0$ such that
\begin{eqnarray}\label{JapanCorrection01}
\epsilon+\delta''\in \left(0,\min\left\{\delta,\int_0^T e^{-\bar{q}_x s} ds \Xi\right\}\right).
\end{eqnarray}
 There is some $\rho\in{\cal R}$ such that
\begin{eqnarray*}\int_0^{T} e^{-{\int_0^t (q_x(\rho_s)-c^G(x,\rho_s))ds }} \int_{\textbf{X}} V^\ast(y)\tilde{q}(dy|x,\rho_t)dt<\infty,~-\infty<\int_0^{T} (q_x(\rho_s)-c^G(x,\rho_s))ds<\infty
\end{eqnarray*}
and
\begin{eqnarray*}
V(x)+\epsilon+\delta''&\ge&\int_0^{T} e^{-\int_0^s (q_x(\rho_v)-c^G(x,\rho_v))dv}\int_{\textbf{X}}V^\ast(y)\tilde{q}(dy|x,\rho_s)ds+e^{-\int_0^T (q_x(\rho_s)-c^G(x,\rho_s))ds}V^\ast(x).
\end{eqnarray*}
Then
\begin{eqnarray*}
\epsilon+\delta''&\ge&\int_0^{T} e^{-\int_0^s (q_x(\rho_v)-c^G(x,\rho_v))dv}\int_{\textbf{X}}V^\ast(y)\tilde{q}(dy|x,\rho_s)ds+e^{-\int_0^T (q_x(\rho_s)-c^G(x,\rho_s))ds}V^\ast(x)-V^\ast(x)\\
&=&\int_0^{T} e^{-\int_0^s (q_x(\rho_v)-c^G(x,\rho_v))dv}\int_{\textbf{X}}V^\ast(y)\tilde{q}(dy|x,\rho_s)ds\\
&&-\int_0^{T} (q_x(\rho_\tau)-c^G(x,\rho_\tau))e^{-\int_0^\tau (q_x(\rho_s)-c^G(x,\rho_s))ds}d\tau V^\ast(x)\\
&=& \int_0^{T} e^{-\int_0^s (q_x(\rho_v)-c^G(x,\rho_v))dv}\left\{\int_{\textbf{X}}V^\ast(y)\tilde{q}(dy|x,\rho_s) - (q_x(\rho_s)-c^G(x,\rho_s))V^\ast(x)\right\}ds\\
&\ge &\int_0^{T} e^{-\int_0^s (q_x(\rho_v)-c^G(x,\rho_v))dv}ds \inf_{a\in\textbf{A}^G}\left\{\int_{\textbf{X}}V^\ast(y)\tilde{q}(dy|x,a) - (q_x(a)-c^G(x,a))V^\ast(x)\right\}\\
&\ge&  \int_0^{T} e^{-\overline{q}_x s}ds \Xi,
\end{eqnarray*}
which is a desired contradiction against (\ref{JapanCorrection01}), meaning that
\begin{eqnarray*}
\Xi=\inf_{a\in\textbf{A}^G}\left\{\int_{\textbf{X}} V^\ast(y) \tilde{q}(dy|x,a)-(q_x(a)-c^G(x,a))V^\ast(x)  \right\}=0.
\end{eqnarray*}}
 $\hfill\Box$
\bigskip

\begin{lemma}\label{JapanLemma05}
Let $w$ be a measurable $[1,\infty)$-valued function satisfying the inequality in Condition \ref{JapanCon02}, whose existence is guaranteed as mentioned in the paragraph below Condition \ref{JapanCon02}. Consider the transition probability $\tilde{p}(dy|x,a)$ on ${\cal B}(\textbf{X})$ given $(x,a)\in\textbf{X}\times\textbf{A}^G$ defined by
\begin{eqnarray*}
\tilde{p}(\Gamma|x,a):=\frac{q(\Gamma|x,a)}{w(x)}+\delta_{x}(dy),~\forall~\Gamma\in{\cal B}(\textbf{X}),~(x,a)\in\textbf{X}\times\textbf{A}^G.
\end{eqnarray*}
Then a $[1,\infty]$-valued lower semianalytic function $V^\ast$ (here the notation $V^\ast$ does not necessarily mean the value function) satisfies (\ref{JapanEqnA}) and (\ref{JapanEqnB}), and for each $x\in\textbf{X}$, either (\ref{JapanEqnA}) or (\ref{JapanEqnB}) holds with equality, if and only if this $[1,\infty]$-valued lower semianalytic function satisfies (\ref{JapanEqnB}), for each $x\in\textbf{X}$
\begin{eqnarray}\label{JapanEqnD}
V^\ast(x)\le \inf_{a\in\textbf{A}^G} \left\{\frac{w(x)}{w(x)-c^G(x,a)}\int_{\textbf{X}}V^\ast(y)\tilde{p}(dy|x,a)\right\},
\end{eqnarray}
and either (\ref{JapanEqnB}) or (\ref{JapanEqnD}) holds with equality, i.e.,
\begin{eqnarray}\label{JapanEqnStar}
V^\ast(x)=\min\left\{\inf_{a\in\textbf{A}^G} \left\{\frac{w(x)}{w(x)-c^G(x,a)}\int_{\textbf{X}}V^\ast(y)\tilde{p}(dy|x,a)\right\},\inf_{b\in\textbf{A}^I}\left\{ \int_{\textbf{X}}V^\ast(y)e^{c^I(x,b,y)}Q(dy|x,b)\right\} \right\}.
\end{eqnarray}
\end{lemma}
Note that (\ref{JapanEqnD}) automatically holds with equality at $x\in \textbf{X}\setminus\textbf{X}^\ast(V^\ast):=\{x\in\textbf{X}:~V^\ast(x)=\infty\}.$ Also note that the function $w$ in the previous lemma does not need be continuous.
\bigskip

\par\noindent\textit{Proof of Lemma \ref{JapanLemma05}.}
``Only if'' part. Consider a $[1,\infty]$-valued lower semianalytic function $V^\ast$ satisfies (\ref{JapanEqnA}) and (\ref{JapanEqnB}), and for each $x\in\textbf{X}$, either (\ref{JapanEqnA}) or (\ref{JapanEqnB}) holds with equality. For $x\in \textbf{X}^\ast(V^\ast)=\{x\in\textbf{X}:~V^\ast(x)<\infty\}$, (\ref{JapanEqnA}) implies for each $a\in\textbf{A}^G$ that
$0\le c^G(x,a)V^\ast(x)+\int_{\textbf{X}}V^\ast(y)q(dy|x,a)=(c^G(x,a)-w(x))V^\ast(x)+w(x)\int_{\textbf{X}}V^\ast(y)\tilde{p}(dy|x,a)$, and thus \begin{eqnarray*}
V^\ast(x)\le \inf_{a\in\textbf{A}^G}\left\{\frac{w(x)}{w(x)-c^G(x,a)}\int_{\textbf{X}}V^\ast(y)\tilde{p}(dy|x,a)\right\},
 \end{eqnarray*}
i.e., (\ref{JapanEqnD}) holds. Let $x\in \textbf{X}^\ast(V^\ast)$ be a point where (\ref{JapanEqnA}) holds with equality. Let us verify at this point $x\in\textbf{X}^\ast(V^\ast),$ (\ref{JapanEqnD}) also holds with equality. For each $\epsilon>0$, there is some $a_\epsilon\in\textbf{A}^G$ such that $\epsilon\ge c^G(x,a_{\epsilon})V^\ast(x)+\int_{\textbf{X}}V^\ast(y)q(dy|x,a_\epsilon)$ so that
\begin{eqnarray*}
V^\ast(x)+\epsilon&\ge& V^\ast(x)+\frac{\epsilon}{w(x)-c^G(x,a_\epsilon)}\ge V^\ast(x)+\frac{c^G(x,a_\epsilon)V^\ast(x)+\int_{\textbf{X}}V^\ast(y)q(dy|x,a_\epsilon)}{w(x)-c^G(x,a_{\epsilon})}\\
&=&\frac{w(x)}{w(x)-c^G(x,a_{\epsilon})}\int_{\textbf{X}}\tilde{p}(dy|x,a_{\epsilon})V^\ast(y)\ge \inf_{a\in\textbf{A}^G} \left\{\frac{w(x)}{w(x)-c^G(x,a)}\int_{\textbf{X}}V^\ast(y)\tilde{p}(dy|x,a)\right\},
\end{eqnarray*}
and thus $V^\ast(x)\ge \inf_{a\in\textbf{A}^G} \left\{\frac{w(x)}{w(x)-c^G(x,a)}\int_{\textbf{X}}V^\ast(y)\tilde{p}(dy|x,a)\right\}.$ The opposite direction of this inequality was seen earlier, and so (\ref{JapanEqnD}) holds with equality at this point. This completes the ``Only if'' part. The argument for the ``If'' part is the same, and omitted. $\hfill\Box$
\bigskip

\begin{remark}\label{JapanRevisionAARemark02}
Consider the function $V^\ast$ in the previous statement. By inspecting the above proof we see the following useful fact: a pair of measurable mappings $\psi^\ast$ and $f^\ast$ from $\textbf{X}$ to $\textbf{A}^I$ and $\textbf{A}^G$ satisfy \begin{eqnarray*}
&&\frac{w(x)}{w(x)-c^G(x,f^\ast(x))}\int_{\textbf{X}}V^\ast(y)\tilde{p}(dy|x,f^\ast(x))=\inf_{a\in\textbf{A}^G} \left\{\frac{w(x)}{w(x)-c^G(x,a)}\int_{\textbf{X}}V^\ast(y)\tilde{p}(dy|x,a)\right\}
\end{eqnarray*}
for each $x\in\textbf{X}$, at which (\ref{JapanEqnD}) holds with equality, and
\begin{eqnarray*} \int_{\textbf{X}} e^{c^I(x,\psi^\ast(x),y)} V^\ast(y)Q(dy|x,\psi^\ast(x)) =\inf_{b\in\textbf{A}^I}\left\{ \int_{\textbf{X}} e^{c^I(x,b,y)}V^\ast(y)Q(dy|x,b) \right\},~\forall~x\in\textbf{X},
\end{eqnarray*}
if and only if
\begin{eqnarray*}
&&\inf_{a\in \textbf{A}^G}\left\{\int_{\textbf{X}}V^\ast(y)\tilde{q}(dy|x,a)-(q_x(a)-c^G(x,a))V^\ast(x)\right\}\\
&=&\int_{\textbf{X}}V^\ast(y)\tilde{q}(dy|x,f^\ast(x))-(q_x(f^\ast(x))-c^G(x,f^\ast(x)))V^\ast(x)
\end{eqnarray*}
for each $x\in\textbf{X}$, at which $0$ coincides with the left hand side, and
\begin{eqnarray*} \int_{\textbf{X}} e^{c^I(x,\psi^\ast(x),y)} V^\ast(y)Q(dy|x,\psi^\ast(x)) =\inf_{b\in\textbf{A}^I}\left\{ \int_{\textbf{X}} e^{c^I(x,b,y)}V^\ast(y)Q(dy|x,b) \right\},~\forall~x\in\textbf{X}.
\end{eqnarray*}
\end{remark}

\begin{lemma}\label{JapanRevisionAAnEWlEMMA}
Conditions \ref{JapanCon02} and \ref{JapanCon01} are satisfied. Then $W^\ast(x)=V^\ast(x)$ for each $x\in\textbf{X}$.
\end{lemma}
\par\noindent\textit{Proof.} According to Proposition \ref{GGZyExponentialProposition02}(a,b), the value function $W^\ast$ for the tilde model is the minimal $[1,\infty]$-valued lower semianalytic function satisfying (\ref{JapanBellman05}) as well as the inequality obtained by replacing the equality in (\ref{JapanBellman05}) by ``$\ge$''. Let us verify that $W^\ast=V^\ast$ as follows.
According to Lemmas \ref{JapanLemma07}, \ref{JapanLemma06} and \ref{JapanLemma05}, the value function $V^\ast$ is a $[1,\infty]$-valued lower semianalytic function satisfying (\ref{JapanBellman05}), c.f. (\ref{JapanEqnStar}). Therefore, $W^\ast\le V^\ast$ pointwise.

For the opposite direction of this inequality, let $x\in\textbf{X}$ be fixed. It suffices to show that $W^\ast$ satisfies (\ref{JapanBellman02}) at the point $x$. Then, since the point $x\in\textbf{X}$ was arbitrarily fixed, one could apply Lemma \ref{JapanNewLem01} to obtain $V^\ast\le W^\ast$ pointwise.

Recall that, as observed in the beginning of this proof, $W^\ast$ satisfies (\ref{JapanEqnStar}). By Lemma \ref{JapanLemma05}, it satisfies (\ref{JapanEqnA}) and (\ref{JapanEqnB}), one of which holds with equality at this point $x.$ If (\ref{JapanEqnB}) holds with equality for $W^\ast$ at $x$, then
\begin{eqnarray*}
&&W^\ast(x)=\inf_{b\in\textbf{A}^I}\left\{\int_{\textbf{X}}W^\ast(y)e^{c^I(x,b,y)}Q(dy|x,b)\right\}\\
&\ge& \inf_{\hat{a}\in\hat{\textbf{A}}}\left\{\int_0^c \int_{\textbf{X}}W^\ast(y)\tilde{q}(dy|x,\rho_t)e^{-\int_0^t (q_x(\rho_s)-c^G(x,\rho_s))ds}dt \right.\\
&&\left.+I\{c=\infty\} e^{-\int_0^\infty q_x(\rho_s)ds} e^{\int_0^\infty c^G(x,\rho_s)ds}\right.\nonumber\\
&&\left.+I\{c<\infty\}e^{-\int_0^c (q_x(\rho_s)-c^G(x,\rho_s))ds}\int_{\textbf{X}}W^\ast(y)e^{c^I(x,b,y)}Q(dy|x,b)\right\},
\end{eqnarray*}
and thus (\ref{JapanBellman02}) is satisfied by $W^\ast$ at $x$, as required. Now suppose (\ref{JapanEqnA}) holds with equality for $W^\ast$ at $x$. It suffices to consider $W^\ast(x)<\infty,$ for otherwise, (\ref{JapanBellman02}) automatically holds for $W^\ast$ at $x.$  According to Remark \ref{JapanRevisionAARemark02} after Lemma \ref{JapanLemma05} and because the tilde model is semicontinuous, there is some $a^\ast\in\textbf{A}^G$ satisfying \begin{eqnarray*}
&&\int_{\textbf{X}} W^\ast(y) \tilde{q}(dy|x,a^\ast)-(q_x(a^\ast)-c^G(x,a^\ast))W^\ast(x)\\
&=&\inf_{a\in\textbf{A}^G}\left\{\int_{\textbf{X}} W^\ast(y) \tilde{q}(dy|x,a)-(q_x(a)-c^G(x,a))W^\ast(x)  \right\}=0,
\end{eqnarray*}
and hence $\int_{\textbf{X}} W^\ast(y) \tilde{q}(dy|x,a^\ast)=(q_x(a^\ast)-c^G(x,a^\ast))W^\ast(x)$. This implies $q_x(a^\ast)\ge c^G(x,a^\ast)$ as the left hand side of the previous equality is nonnegative and $W^\ast(x)\ge 1,$ and for the same reason, if $c^G(x,a^\ast)=q_x(a^\ast)$, then $c^G(x,a^\ast)=q_x(a^\ast)=0$, in which case,
 \begin{eqnarray*}
&&W^\ast(x)\ge 1=  \int_0^\infty \int_{\textbf{X}}W^\ast(y)\tilde{q}(dy|x,a^\ast)e^{-\int_0^t (q_x(a^\ast)-c^G(x,a^\ast))ds}dt+ e^{-\int_0^\infty q_x(a^\ast)ds} e^{\int_0^\infty c^G(x,a^\ast)ds}\\
&\ge& \inf_{\hat{a}\in\hat{\textbf{A}}}\left\{\int_0^c \int_{\textbf{X}}W^\ast(y)\tilde{q}(dy|x,\rho_t)e^{-\int_0^t (q_x(\rho_s)-c^G(x,\rho_s))ds}dt \right.\\
&&\left.+I\{c=\infty\} e^{-\int_0^\infty q_x(\rho_s)ds} e^{\int_0^\infty c^G(x,\rho_s)ds}\right.\nonumber\\
&&\left.+I\{c<\infty\}e^{-\int_0^c (q_x(\rho_s)-c^G(x,\rho_s))ds}\int_{\textbf{X}}W^\ast(y)e^{c^I(x,b,y)}Q(dy|x,b)\right\}.
\end{eqnarray*}
That is, (\ref{JapanBellman02}) is satisfied by $W^\ast$ at $x$, as desired. Finally, if $c^G(x,a^\ast)<q_x(a^\ast)$, then
\begin{eqnarray*}
 &&\inf_{\hat{a}\in\hat{\textbf{A}}}\left\{\int_0^c \int_{\textbf{X}}W^\ast(y)\tilde{q}(dy|x,\rho_t)e^{-\int_0^t (q_x(\rho_s)-c^G(x,\rho_s))ds}dt \right.\\
&&\left.+I\{c=\infty\} e^{-\int_0^\infty q_x(\rho_s)ds} e^{\int_0^\infty c^G(x,\rho_s)ds}\right.\nonumber\\
&&\left.+I\{c<\infty\}e^{-\int_0^c (q_x(\rho_s)-c^G(x,\rho_s))ds}\int_{\textbf{X}}W^\ast(y)e^{c^I(x,b,y)}Q(dy|x,b)\right\}\\ &\le&\int_0^\infty \int_{\textbf{X}}W^\ast(y)\tilde{q}(dy|x,a^\ast)e^{-\int_0^t (q_x(a^\ast)-c^G(x,a^\ast))ds}dt+ e^{-\int_0^\infty q_x(a^\ast)ds} e^{\int_0^\infty c^G(x,a^\ast)ds}\\
 &=&\frac{\int_{\textbf{X}}W^\ast(y)\tilde{q}(dy|x,a^\ast)}{q_x(a^\ast)-c^G(x,a^\ast)}+0=W^\ast(x),
\end{eqnarray*}
as requested. Thus, $W^\ast$ satisfies (\ref{JapanBellman02}). Consequently, $W^\ast= V^\ast$ on $\textbf{X}$, as required. $\hfill\Box$
\bigskip

\par\noindent\textit{Proof of Theorem \ref{JapanTheorem05}}.  Part (b) was seen in the proof of Lemma \ref{JapanLemma06}.

Consider the pair of measurable mappings $(\psi^\ast,f^\ast)$ from Proposition \ref{JapanRevisionAAProp01}. Recall that $W^\ast=V^\ast$ on $\textbf{X}$ by Lemma \ref{JapanRevisionAAnEWlEMMA}. Keeping in mind Remark \ref{JapanRevisionAARemark02}, an inspection of the proof of Lemma \ref{JapanRevisionAAnEWlEMMA} reveals that the deterministic stationary strategy $(\varphi(x),\psi^\ast(x),t\rightarrow \delta_{f^\ast(x)}(da))\in\hat{\textbf{A}}$ in the hat DTMDP model, where $\varphi$ is defined in part (c) of this theorem, attains the infimum in
 \begin{eqnarray*}
V^\ast(x)&=&\inf_{\hat{a}\in\hat{\textbf{A}}}\left\{\int_0^c \int_{\textbf{X}}V^\ast(y)\tilde{q}(dy|x,\rho_t)e^{-\int_0^t (q_x(\rho_s)-c^G(x,\rho_s))ds}dt \right.\\
&&\left.+I\{c=\infty\} e^{-\int_0^\infty q_x(\rho_s)ds} e^{\int_0^\infty c^G(x,\rho_s)ds}\right.\nonumber\\
&&\left.+I\{c<\infty\}e^{-\int_0^c (q_x(\rho_s)-c^G(x,\rho_s))ds}\int_{\textbf{X}}V^\ast(y)e^{c^I(x,b,y)}Q(dy|x,b)\right\}\nonumber
\end{eqnarray*}
for each $x\in\textbf{X}$. By Theorem \ref{JapanTheorem01}, this deterministic stationary strategy $(\varphi(x),\psi^\ast(x),t\rightarrow \delta_{f^\ast(x)}(da))\in\hat{\textbf{A}}$ is optimal for problem (\ref{JapanProblem02}) for the hat DTMDP model. This and Remark \ref{JapanRevisionAARemark01} imply that $V^\ast={\cal V}^\ast$ on $\textbf{X}$ and part (c). By Lemma \ref{JapanRevisionAAnEWlEMMA}, we see now ${\cal V}^\ast=W^\ast$ on $\textbf{X}$, and thus part (a) holds.
$\hfill\Box$
\bigskip

\par\noindent\textit{Proof of Corollary \ref{JapanRevisionCorollary0001}.} This corollary follows at once from Theorem  \ref{JapanTheorem05}, Lemma \ref{JapanLemma05} and Remark \ref{JapanRevisionAARemark02}. $\hfill\Box$
\bigskip

% and the first assertion of Theorem \ref{JapanTheorem05} is seen. The last assertion of Theorem \ref{JapanTheorem05} follows from the paragraph after Lemma \ref{JapanLemma05}, as noted earlier in this proof.

%Now we turn to the proof of Theorem \ref{JapanTheorem01}. Since $W^\ast$ is lower semicontinuous, so is $V^\ast=W^\ast$, according to Theorem \ref{JapanTheorem05}, just verified. Part (a) of Theorem \ref{JapanTheorem01} is seen.  %It follows from the semicontinuity of the tilde model the existence of a measurable selectors $(\psi^\ast,f^\ast)$ satisfying (\ref{JapanConserving01}).
%This pair of measurable selectors satisfies the relations in part (c) of Theorem \ref{JapanTheorem01} according to Theorem \ref{JapanTheorem05}. The deterministic stationary strategy in part (d) is optimal because of  Lemma \ref{JapanNewLem01} and that This fact can be seen. Thus, Theorem \ref{JapanTheorem01} is proved. $\hfill\Box$
%\bigskip

\appendix
\section{Appendix: relevant results about DTMDPs}\label{JapanAppendixA}

In this appendix, we present the relevant facts about DTMDPs. The proofs of the presented statements can be found in \cite{Jaskiewicz:2008} or \cite{Zhang:2017}. Standard description of a DTMDP can be found in e.g., \cite{Hernandez-Lerma:1996,Piunovskiy:1997}. The notations used in this section are independent of the previous sections.

A DTMDP has the following primitives $\{\textbf{X},\textbf{A}, p, $l$\}$:
\begin{itemize}
\item $\textbf{X}$ is a nonempty Borel state space.
\item $\textbf{A}$ is a nonempty Borel action space.
\item $p(dy|x,a)$ is a stochastic kernel on ${\cal B}(\textbf{X})$ given $(x,a)\in \textbf{X}\times\textbf{A}$.
\item $l$ a $[0,\infty]$-valued measurable cost function on $\textbf{X}\times\textbf{A}\times\textbf{X}.$
\end{itemize}

\begin{condition}\label{GGZyExponentialDTMDPCon2}
\begin{itemize}
\item[(a)] The function $l(x,a,y)$ is lower semicontinuous in  $(x,a,y)\in \textbf{X}\times\textbf{A}\times\textbf{X}.$
\item[(b)] For each bounded continuous function $f$ on $\textbf{X}$, $\int_{\textbf{X}}f(y)p(dy|x,a)$ is continuous in $(x,a)\in\textbf{X}\times \textbf{A}.$
\item[(c)] The space $\textbf{A}$ is a compact Borel space.
\end{itemize}
\end{condition}

\begin{definition}\label{JapanRevisionAADef01}
The DTMDP model $\{\textbf{X},\textbf{A}, p, $l$\}$ is called semicontinuous if it satisfies Condition \ref{GGZyExponentialDTMDPCon2}.
\end{definition}

Let us denote for each $n=1,2,\dots,\infty,$ $\textbf{H}_n:=\textbf{X}\times(\textbf{A}\times\textbf{X})^n$ and $\textbf{H}_{0}:=\textbf{X}.$
A strategy $\sigma=(\sigma_n)_{n=0}^\infty$ in the DTMDP is given by a sequence of stochastic kernels $\sigma_n(da|h_{n})$ on ${\cal B}(\textbf{A})$ from $h_{n}\in \textbf{H}_{n}$ for $n=0,1,2,\dots.$ A strategy $\sigma=(\sigma_n)$ is called deterministic Markov if for each $n=0,1,2,\dots,$ $\sigma_n(da|h_{n})=\delta_{\{\varphi_n(x_{n})\}}(da)$, where $\varphi_{n}$ is an $\textbf{A}$-valued measurable mapping on $\textbf{X}.$ We identify such a deterministic Markov strategy with $(\varphi_n).$ A deterministic Markov strategy $(\varphi_n)$ is called deterministic stationary if $\varphi_n$ does not depend on $n,$ and it is identified with the underlying measurable mapping $\varphi$ from $\textbf{X}$ to $\textbf{A}.$ Let $\Sigma$ be the space of strategies, and $\Sigma_{DM}$ be the space of all deterministic strategies for the DTMDP.

Let the controlled and controlling process be denoted by $\{Y_n, n=0,1,\dots,\infty\}$ and $\{A_n,n=0,1,\dots,\infty\}$. Here, for each $n=0,1,\dots,$ $Y_n$ is the projection of $\textbf{H}_\infty$ to the $2n+1$st coordinate, and $A_n$ to the $2n+2$nd coordinate.
Under a strategy $\sigma=(\sigma_n)$ and a given initial probability distribution $\nu$ on $(\textbf{X},{\cal B}(\textbf{X}))$, by the Ionescu-Tulcea theorem, c.f., \cite{Hernandez-Lerma:1996,Piunovskiy:1997}, one can construct a probability measure $\textbf{P}_\nu^\sigma$ on $(\textbf{H}_\infty,{\cal B}(\textbf{H}_\infty))$ such that
\begin{eqnarray*}
&&\textbf{P}_\nu^\sigma(Y_0\in dx)=\nu(dx),\\
&&\textbf{P}_\nu^\sigma(A_n\in da|Y_0,A_0,\dots,Y_n)=\sigma_{n}(da|Y_0,A_0,\dots,Y_n),~n=0,1,\dots,\\
&&\textbf{P}_\nu^\sigma(Y_{n+1}\in dx|Y_0,A_0,\dots,Y_n,A_n)=p(dx|Y_n,A_n),~n=0,1,\dots.
\end{eqnarray*}
As usual, equalities involving conditional expectations and probabilities are understood in the almost sure sense.
\begin{definition}\label{JapanRevisionDef01}
The probability measure $\textbf{P}_\nu^\sigma$ is called a strategic measure (of the strategy $\sigma$) in the DTMDP model $\{\textbf{X},\textbf{A}, p, $l$\}$ (with the initial distribution $\nu$).
\end{definition}
The expectation taken with respect to $\textbf{P}_\nu^\sigma$ is denoted by $\textbf{E}_\nu^\sigma.$ When $\nu$ is concentrated on the singleton $\{x\}$, $\textbf{P}_\nu^\sigma$ and $\textbf{E}_\nu^\sigma$ are written as $\textbf{P}_x^\sigma$ and $\textbf{E}_x^\sigma.$

Consider the optimal control problem
\begin{eqnarray}\label{ZyExponentialProblem2}
\mbox{Minimize over $\sigma$}:&& \textbf{E}_x^\sigma\left[e^{\sum_{n=0}^\infty l(Y_n,A_n,Y_{n+1})}\right]=:\textbf{V}(x,\sigma),~x\in \textbf{X}.
\end{eqnarray}
We denote the value function of problem (\ref{ZyExponentialProblem2}) by $\textbf{V}^\ast$. Then a strategy $\sigma^\ast$ is called optimal for problem  (\ref{ZyExponentialProblem2}) if $\textbf{V}(x,\sigma^\ast)=\textbf{V}^\ast(x)$ for each $x\in \textbf{X}.$ For a constant $\epsilon>0$, a strategy is called $\epsilon$-optimal for problem  (\ref{ZyExponentialProblem2}) if $\textbf{V}(x,\sigma^\ast)\le\textbf{V}^\ast(x)+\epsilon$ for each $x\in \textbf{X}.$

Occasionally we will also consider the so called universally measurable strategies, in which case, the stochastic kernels $\sigma_n(da|h_{n})$ are universally measurable, i.e., for each measurable subset $\Gamma$ of $\textbf{A}$, $\sigma(\Gamma|h_n)$ is universally measurable in $h_n\in\textbf{H}_n.$ The meaning of universally measurable deterministic Markov or deterministic stationary strategy is understood similarly, i.e., when the underlying mappings are universally measurable in their arguments. See Chapter 7.7 of \cite{Bertsekas:1978} for the definition of universal measurability and other related measurability concepts, such as the definition of a lower semianalytic function.

We collect the relevant statements in Section 3 of \cite{Zhang:2017} in the next proposition.
\begin{proposition}\label{GGZyExponentialProposition02}
The following assertions hold.
\begin{itemize}
\item[(a)] The value function $\textbf{V}^\ast$ is the minimal $[1,\infty]$-valued lower semianalytic solution to \begin{eqnarray}\label{ZyExponential02}
\textbf{V}(x)=\inf_{a\in \textbf{A}}\left\{\int_{\textbf{X}}e^{l(x,a,y)}\textbf{V}(y)p(dy|x,a)\right\},~x\in \textbf{X}.
\end{eqnarray}
\item[(b)] Let $\textbf{U}$ be a $[1,\infty]$-valued lower semianalytic function on $\textbf{X}$. If
\begin{eqnarray*}
\textbf{U}(x)\ge \inf_{a\in \textbf{A}}\left\{\int_{\textbf{X}}e^{l(x,a,y)}\textbf{U}(y)p(dy|x,a)\right\},~\forall~x\in \textbf{X},
\end{eqnarray*}
then $\textbf{U}(x)\ge \textbf{V}^\ast(x)$ for each $x\in \textbf{X}.$
\item[(c)] Let $\varphi$ be a deterministic stationary strategy for the DTMDP model $\{\textbf{X},\textbf{A},p,l\}$. If
\begin{eqnarray}\label{GGZyExponential025}
\textbf{V}^\ast(x)=\int_{\textbf{X}}e^{l(x,\varphi(x),y)}\textbf{V}^\ast(y)p(dy|x,\varphi(x)),~\forall~x\in \textbf{X},
\end{eqnarray}
then $\textbf{V}^\ast(x)=\textbf{V}(x,\varphi)$ for each $x\in \textbf{X}.$
\item[(d)] $\textbf{V}^\ast(x)=\inf_{\sigma\in \Sigma^U}\textbf{V}(x,\sigma),$ where $\Sigma^U$ is the set of universally measurable strategies. Moreover, for each $\epsilon>0$, there is some universally measurable deterministic stationary $\epsilon$-optimal strategy for problem (\ref{ZyExponentialProblem2}).
\item[(e)] Suppose Condition \ref{GGZyExponentialDTMDPCon2} is satisfied. Then the value function $\textbf{V}^\ast$ is the minimal $[1,\infty]$-valued lower semicontinuous solution to (\ref{ZyExponential02}). Moreover, there exists a deterministic stationary strategy $\varphi$ satisfying (\ref{GGZyExponential025}), and so in particular, there exists a deterministic stationary optimal strategy for problem (\ref{ZyExponentialProblem2}).
\end{itemize}
\end{proposition}
Part (d) of the above statement follows from the proof of Proposition 3.2 of \cite{Zhang:2017}, whereas all the other parts are according to Propositions 3.1, 3.4 and 3.7 therein.

%Stronger results than those of Proposition \ref{GGZyExponentialProposition02} hold for semicontinuous model. The next statement is taken from Proposition 3.7 of \cite{Zhang:2017}.
%\begin{proposition}\label{JapanAppendixProp01}
%Suppose Condition \ref{GGZyExponentialDTMDPCon2} is satisfied. Then the following assertions hold.
%\begin{itemize}
%\item[(a)] The value function $\textbf{V}^\ast$ is the minimal $[1,\infty]$-valued lower semicontinuous solution to (\ref{ZyExponential02}). Moreover, there exists a deterministic stationary strategy $\varphi$ satisfying (\ref{GGZyExponential025}), and so in particular, there exists a deterministic stationary optimal strategy for problem (\ref{ZyExponentialProblem2}).

%\item[(b)] Let $\textbf{V}^{(0)}(x):=1$ for each $x\in \textbf{X}$, and for each $n=1,2,\dots,$
 %   \begin{eqnarray*}
  %  \textbf{V}^{(n)}(x):=\inf_{a\in A}\left\{\int_{\textbf{X}}p(dy|x,a)e^{l(x,a,y)}\textbf{V}^{(n-1)}(y)\right\},~\forall~x\in \textbf{X}.
  %  \end{eqnarray*}
%Then $(\textbf{V}^{(n)}(x))$ increases to $\textbf{V}^\ast(x)$ for each $x\in \textbf{X}$.
%\end{itemize}
%\end{proposition}

\subsection*{Acknowledgement}
We thank the editor and referees for comments and remarks that improve significantly the readability of this paper. This work is supported by the Royal Society (grant number IE160503), and the Daiwa Anglo-Japanese Foundation (UK) (grant reference 4530/12801).

\end{document}